\documentclass[11pt]{article}
\usepackage[T1]{fontenc}
\usepackage[utf8]{inputenc}
\usepackage{amsmath,amssymb,amsthm,microtype}
\usepackage[hidelinks]{hyperref}
\newtheorem{theorem}{Theorem}[section]
\newtheorem{lemma}[theorem]{Lemma}
\newtheorem{proposition}[theorem]{Proposition}
\newtheorem{corollary}[theorem]{Corollary}
\theoremstyle{definition}
\newtheorem{definition}[theorem]{Definition}
\newtheorem{remark}[theorem]{Remark}
\newtheorem{conjecture}[theorem]{Conjecture (unproved)}
\newtheorem{question}[theorem]{Open question}
\newcommand{\R}{\mathbb R}
\newcommand{\Q}{\mathbb Q}
\newcommand{\Z}{\mathbb Z}
\newcommand{\T}{\mathbb T}
\newcommand{\tu}{\tau_{\mathrm{usual}}}
\newcommand{\dist}{\operatorname{dist}}
\title{Almost Periodic Characters and Hamel Coefficients in TVG Topologies on the Real Line}
\author{Irina Yaroshevskaya}
\date{}
\hypersetup{pdftitle={Almost Periodic Characters and Hamel Coefficients in TVG Topologies on the Real Line},pdfauthor={Irina Yaroshevskaya}}
\begin{document}
\maketitle
\begin{abstract}
For a nonzero rational-valued additive function $f$ on $\R$, modular Hamel observations generate topologies $\beta_f$ and $\sigma_f$, while real-valued observations generate $\tau_f$. We prove $\tu\subsetneq\sigma_f\subsetneq\tau_f$, determine the characters and real-valued additive functions of the modular pair, and classify their equality and continuous homomorphisms. Each modular family contains $2^{\mathfrak c}$ pairwise nonisomorphic topological groups. Their completions are $\R_{\mathrm{usual}}\times\widehat{\R_d}$ and $\widehat{\R_d}\times\widehat{\R_d}$, respectively. The topologies $\beta_f$ and $\sigma_f$ have the same compact subsets and convergent sequences. Nevertheless the explicit sequence $n!(3-2\sqrt2)^{n^2}$ converges in both and fails to converge in $\tau_f$ for every nonzero $f$. We classify convergence in finite rational spans and exhibit different compact subsets for $\sigma_f$ and $\tau_f$. A character-subspace description classifies all Hausdorff precompact TVG topologies. Using de la Barrera Mayoral's compatible refinement on $\Q/\Z$ \cite{Barrera2014}, we construct a locally quasi-convex TVG strictly finer than $\sigma_f$ with the same character set. Finally, two Hausdorff nondiscrete TVGs strictly finer than usual are constructed with discrete join. Thus no greatest nondiscrete TVG exists, and maximal nondiscrete extensions are not unique. For the separate greatest-compatible-topology question, we give explicit rational duality tests, give ultrafilter criteria, and construct incomparable compatible rational topologies from prime square denominators. A quotient-lifting argument gives a sufficient counterexample criterion on the torsion group $\Q/\Z$. We apply established polar criteria and prove that two nuclear compatible topologies have a compatible join. The prime-pair join character set, nuclearity of its members, and existence of a greatest compatible locally quasi-convex TVG remain unresolved here. Standard duality and cited background results are distinguished from the specialized arguments; no historical priority claim is made. Real dilations are the only automorphisms of the Hamel dual pair, so a greatest unrestricted compatible topology would automatically be a TVG topology and would settle the Greatest Compatible TVG Conjecture; a real saturation criterion and a necessary condition for a greatest member are proved, and a dichotomy shows that no splitting of the dyadic chain yields an incompatible pair.
\end{abstract}

\clearpage
\tableofcontents
\clearpage

\section{Introduction and previous results}
Topological vector group topologies on the real line can be studied using additive functions with values either in the real line or in the circle group. These two target groups encode different information. Real-valued observations control numerical size, whereas circle-valued observations control distance from integers. This article gives an explicit separation between the resulting topologies, even after adjoining the usual topology.

In \cite{Y1975} the author introduced a Hausdorff quasi-vector topology $t$ on $\R$ different from the usual topology. Its minimality among Hausdorff topological vector group topologies was proved in \cite{Y1981}, received in 1979 and published in 1981. Here ``minimal'' means that no strictly coarser Hausdorff TVG topology exists; it is not used as a synonym for ``least'' in an arbitrary class. The recent article \cite{Y2026} constructs topologies $\rho_{(r_n)}$ such that
\[
t\subsetneq\rho_{(r_n)}\subsetneq\tu
\]
for positive sequences satisfying $\frac{-\log r_n}{n}\to\infty$. Factorial powers, quadratic exponential decay, and prime radicals divided by a quadratic exponential provide examples there.

The companion manuscript \cite{YHamel} treats rational-valued Hamel coefficient topologies above $\tu$. It proves an exact comparison theorem, $2^{\mathfrak c}$ distinct topologies on the fixed real line, failure of first countability, and a common nondiscrete refinement, where $\mathfrak c=|\R|$. That article is available as arXiv:2610.24132; its results are recalled here as previous work. The present article focuses on a different construction: modular Hamel observations produce topologies strictly between $\tu$ and those coefficient topologies, but with no nonzero continuous rational-valued additive functions. We also determine their continuous characters and real-valued additive functions, classify equality and inclusion, identify their group completions, and settle first countability and local compactness.

Bohr topologies, initial topologies, and the use of Hamel functions are classical. Dom\'inguez and Tarieladze \cite{DT2008} discuss real-line Bohr topologies in Example~2.6 and additive versus real-linear character lifts in Section~4. Dergachev \cite{D1978} discusses the almost-periodic real-line example in the TVG setting. Hejcman \cite{H1988} studies TVG topologies on the real line. Our contribution here is presented as a self-contained theorem and proof; historical novelty has not been established by the literature check, and no priority relative to Hejcman's paper is asserted.

The additional results concern six directions: homomorphism rigidity, compatible topologies, compact subsets and sequences, group isomorphisms, precompact classification, and maximal nondiscrete extensions. Sections~\ref{sec:rigidity}--\ref{sec:maximal} give the proved results and state the limits of the investigation. In particular, an explicit sequence separates modular from real-valued observations for every nonzero coefficient function, while two nondiscrete scalar-invariant topologies with discrete join exclude a greatest nondiscrete TVG. A compatible refinement also rules out $\sigma_f$ as the finest topology for its character set; the existence of some other greatest compatible member is not settled. Results using Glicksberg's theorem or standard precompact duality are presented as applications of those classical tools, without a novelty claim.

The compatible-topology investigation is developed further in Sections~\ref{sec:joins-v8}--\ref{sec:status-v8}. Directed compatible joins and maximal compatible extensions are established, rational quotients and transfer maps are analyzed, and ultrafilter tests distinguish algebraic limits from continuous characters. The explicit rational topologies $\gamma_2,\gamma_3$ are compatible and incomparable; their join is weaker than usual, but its full character set has not been identified. The finite-order obstruction and the unsuccessful block construction documented in Appendix~\ref{app:blocks-v11} explain why those examples do not prove nonexistence of a greatest compatible TVG. Appendix~\ref{app:segments-v12} proves discreteness of a factorial-segment family and states the limits of finite numerical tests. The main open existence assertion is named the Greatest Compatible TVG Conjecture (Conjecture~\ref{conj:greatest-v12}). Sections~\ref{sec:scalar-v10}--\ref{sec:polar-v10} distinguish scalar invariance from compatibility, lift compatible topologies from $\Q/\Z$, and give the precise established polar and nuclearity criteria. All counterexample implications there remain conditional on constructing an incompatible pair. Section~\ref{sec:dilations-v13} shows that real dilations are the only automorphisms of the dual pair, so that a greatest unrestricted compatible topology would automatically be a TVG topology; this makes the Greatest Compatible TVG Conjecture a consequence of the Mackey problem for the Hamel pair, and it records a real saturation criterion and a necessary condition. Appendix~\ref{app:dyadic-v13} proves that no splitting of the dyadic chain gives an incompatible pair.

\section{Definitions and notation}
Write $\Z$ for the integers, $\Q=\{p/q:p,q\in\Z,\ q\ne0\}$ for the rational numbers, and $\R$ for the real numbers. The usual topology on $\R$ is $\tu$. A topology is stronger, or finer, than another if it contains all its open sets. The join $\tau\vee\eta$ is the smallest topology containing both $\tau$ and $\eta$.

\begin{definition}
A \emph{topological vector group (TVG) topology} on $\R$ is a group topology for addition such that $x\mapsto ax$ is continuous for every fixed $a\in\R$. Joint continuity of $(a,x)\mapsto ax$, with the scalar variable carrying its usual topology, is not required. Unless stated otherwise, the TVG topologies considered here are Hausdorff.
\end{definition}
For $a\ne0$ the dilation is a homeomorphism, since its inverse is multiplication by $1/a$. This convention agrees with \cite[Definition~2.1]{DT2008}.

\begin{definition}
A \emph{Hamel basis} $H$ of $\R$ over $\Q$ is a rationally linearly independent spanning set. Every real number has a unique expression
\[
x=\sum_{h\in F}q_hh,\qquad F\subset H\text{ finite},\quad q_h\in\Q,
\]
where omitted coefficients are zero. Choosing $h_0\in H$ defines a \emph{coefficient function} $f(x)=q_{h_0}$. It is additive and rational-linear.
\end{definition}
Existence of Hamel bases uses the axiom of choice. Their connection with discontinuous additive functions goes back to Hamel \cite{Hamel1905}.
\[
D=\operatorname{Hom}_{\Q}(\R,\Q).
\]
Here $\operatorname{Hom}_{\Q}$ denotes rational-linear maps. Every additive map from $\R$ to $\Q$ is rational-linear. Every nonzero $f\in D$ is a coefficient function for some Hamel basis: choose $h_0$ with $f(h_0)=1$, take a basis of $\ker f$, and adjoin $h_0$. The kernel is
\[
\ker f=\{x\in\R:f(x)=0\}.
\]
It equals the rational span of the unselected basis elements in this representation. We use the notation
\[
f_a(x)=f(ax)\qquad(a\in\R).
\]
The usual topology on $\Q$ is its subspace topology inherited from $\R$.

Let $\T=\R/\Z$ carry its usual circle topology, and write
\[
p(r)=r+\Z,\qquad \|r\|_{\T}=\dist(r,\Z)=\inf_{m\in\Z}|r-m|.
\]
A \emph{character} on an additive group is a homomorphism into $\T$. In particular $p\circ g$ is a character for every additive real-valued $g$.

The \emph{initial topology} of a family of maps into specified topological spaces is the coarsest topology making every map continuous. Finite intersections of inverse images of open sets form a base. Our maps are additive, so their initial topologies are group topologies. We write
\[
D_{\tau}=\operatorname{Hom}_{\mathrm{cont}}((\R,\tau),\Q_{\mathrm{usual}})
\]
for the set of continuous additive maps to the usual rational numbers. The equation $D_\tau=\{0\}$ means that only the zero function is continuous. This is not the real-linear dual and is not the group of continuous characters.

\section{Almost periodic observations and the Bohr topology}
A usual-continuous complex-valued function on $\R$ is \emph{Bohr uniformly almost periodic} if it is a uniform limit of trigonometric polynomials
\[
P(x)=\sum_{j=1}^m c_j e^{2\pi i b_jx},\qquad c_j\in\mathbb C,\quad b_j\in\R.
\]
This standard description makes its connection with characters explicit. Let $t$ be the initial topology of the characters
\[
x\longmapsto bx+\Z\qquad(b\in\R).
\]
Its zero-neighborhoods have a base
\[
T(B,\varepsilon)=\{x:\|bx\|_{\T}<\varepsilon\text{ for all }b\in B\},
\]
with $B\subset\R$ finite and $\varepsilon>0$. The topology generated by all usual-continuous uniformly almost periodic functions is the same topology: characters are among those functions, and their uniform limits remain continuous for the character topology. This is the Bohr topology of the usual additive real line, also described in \cite[Example~2.6]{DT2008} and \cite{Y2026}.

A group topology is \emph{precompact} if, for every zero-neighborhood $U$, finitely many translates of $U$ cover the group. The character map embeds $(\R,t)$ into a product of compact circles, so $t$ is precompact. Ordinary characters separate points, making $t$ Hausdorff. The scalar identity $p(bcx)=p((bc)x)$ proves the TVG property. Since the usual real line is not precompact, $t\subsetneq\tu$. As a continuous image of the usual real line, $(\R,t)$ is path connected.

\section{Real-valued and modular Hamel topologies}
Fix $f\in D\setminus\{0\}$. Define $\tau_f$ as the initial topology of the identity into $\R_{\mathrm{usual}}$ and all $f_a$ into $\Q_{\mathrm{usual}}$. A zero-neighborhood base is
\[
U_f(\delta,A,\varepsilon)
=\{x:|x|<\delta,\ |f(ax)|<\varepsilon\ (a\in A)\}.
\]
Here $\delta,\varepsilon>0$ and $A\subset\R$ is finite, possibly empty.

Define $\beta_f$ as the initial topology of all ordinary characters $p(bx)$ and all Hamel characters $p(f(ax))$, for $a,b\in\R$. Explicitly, a base at zero consists of
\[
W_f(A,B,\varepsilon)
=\{x:\|f(ax)\|_{\T}<\varepsilon\ (a\in A),\
\|bx\|_{\T}<\varepsilon\ (b\in B)\},
\]
where $A,B\subset\R$ are finite, possibly empty, and $\varepsilon>0$. Thus $\beta_f$ controls residues modulo integers rather than the real sizes of its observations. Finally, set
\[
\sigma_f=\tu\vee\beta_f.
\]
Since the ordinary characters are already usual-continuous, a base at zero for $\sigma_f$ consists of
\begin{equation}\label{eq:base}
V_f(\delta,A,\varepsilon)
=\{x:|x|<\delta,\ \|f(ax)\|_{\T}<\varepsilon\ (a\in A)\},
\end{equation}
where $\delta,\varepsilon>0$ and $A$ is finite, possibly empty. A common positive bound suffices by taking the minimum of finitely many bounds. Neighborhoods at any $y$ are translates of these zero-neighborhoods. In the base for $\tau_f$, the inequalities $|f(ax)|<\varepsilon$ use the usual absolute value on rational numbers; in the bases for $\beta_f$ and $\sigma_f$, $\|f(ax)\|_{\T}<\varepsilon$ instead permits values close to any integer. This distinction is responsible for the strict inclusion proved below.

\begin{lemma}\label{lem:noncontinuous}
For nonzero $g\in D$, neither $g$ nor $p\circ g$ is usual-continuous.
\end{lemma}
\begin{proof}
A usual-continuous additive real function has the form $x\mapsto cx$. Its image can lie in $\Q$ only if $c=0$, proving the first assertion. Every usual-continuous character on $\R$ has the form $x\mapsto p(cx)$; see \cite[Proposition~1.3]{DT2008}. If $p\circ g=p(c\,\cdot)$, the additive map $g(x)-cx$ takes integer values. Every additive map $h:\R\to\Z$ is zero: for each positive integer $n$, $h(x)=nh(x/n)$, so $h(x)$ is divisible by every $n$. Hence $g(x)=cx$, a contradiction.
\end{proof}

\begin{proposition}\label{prop:basic}
The topologies $\beta_f$, $\sigma_f$, and $\tau_f$ are Hausdorff nondiscrete TVG topologies. The topology $\beta_f$ is precompact, and
\[
t\subsetneq\beta_f,\qquad
\tu\subsetneq\sigma_f\subseteq\tau_f.
\]
The topologies $\beta_f$ and $\tu$ are incomparable. The strict inclusion $\sigma_f\subsetneq\tau_f$ will be established separately in Theorem~\ref{thm:main}.
\end{proposition}
\begin{proof}
The generating maps are homomorphisms, so all three topologies are group topologies. Composition with multiplication by $c$ replaces each parameter $a$ or $b$ by $ac$ or $bc$; the identity coordinate also remains continuous. This proves the TVG property. Ordinary characters separate points in $\beta_f$; the other two topologies contain $\tu$, so all are Hausdorff.

For every fixed $a$, $f(a/n)=f(a)/n\to0$. Together with $1/n\to0$ usually, this proves that $1/n\to0$ in $\tau_f$, and therefore also in $\sigma_f$ and $\beta_f$. None is discrete. The character embedding places $(\R,\beta_f)$ in a product of compact circles, proving precompactness.

All defining characters for $\beta_f$ are $\tau_f$-continuous, so $\sigma_f\subseteq\tau_f$. Lemma~\ref{lem:noncontinuous} shows that $p\circ f$ is continuous in $\beta_f$ but not in $\tu$. Thus $\beta_f\not\subseteq\tu$, $t\subsetneq\beta_f$, and $\tu\subsetneq\sigma_f$. Conversely, $\tu\subseteq\beta_f$ would make the usual real line a continuous homomorphic image of a precompact group and hence precompact, which is false. This proves incomparability.
\end{proof}

\section{The finite-kernel argument}
The main proof applies more generally than the single-function construction. Let $\mathcal S\subseteq D$ be any family, and let $\sigma_{\mathcal S}$ be generated by $\tu$ and the characters $p\circ s$, for $s\in\mathcal S$. A base is
\[
V(\delta;s_1,\ldots,s_m;\varepsilon)
=\{x:|x|<\delta,\ \|s_j(x)\|_{\T}<\varepsilon\ (1\le j\le m)\}.
\]
No invariance under scalars is needed for the next theorem.

\begin{theorem}[Unboundedness and vanishing rational-valued dual]\label{thm:vanishing}
For every family $\mathcal S\subseteq D$ and every nonzero $g\in D$, the function $g$ is unbounded in absolute value on every $\sigma_{\mathcal S}$-neighborhood of zero. Consequently
\[
D_{\sigma_{\mathcal S}}=\{0\}.
\]
\end{theorem}
\begin{proof}
It suffices to consider a basic neighborhood $V$ as displayed above. Put
\[
K=\ker g\cap\bigcap_{j=1}^m\ker s_j.
\]
The rational-linear map $x\mapsto(g(x),s_1(x),\ldots,s_m(x))$ has values in the finite-dimensional space $\Q^{m+1}$. Therefore $K$ has finite codimension over $\Q$ in $\R$. Since $\R$ has infinite rational dimension, $K\ne\{0\}$. Choose $k\in K\setminus\{0\}$. The set $k\Q$ is usual-dense in $\R$ and lies in $K$, so $K$ is usual-dense.

Choose $y\in\R$ with $g(y)\ne0$. The numbers $s_j(y)$ are rational; choose a positive integer $d$ clearing their finitely many denominators. For every positive integer $N$, density of $K$ supplies $k_N\in K$ for which
\[
x_N=Ndy+k_N,\qquad |x_N|<\delta.
\]
For each $j$ we have
\[
s_j(x_N)=Nd\,s_j(y)\in\Z.
\]
Thus $x_N\in V$, with all modular observations exactly zero. However,
\[
|g(x_N)|=Nd\,|g(y)|\longrightarrow\infty.
\]
This proves unboundedness on $V$, and hence on every zero-neighborhood. If $g$ were continuous into $\Q_{\mathrm{usual}}$, the inverse image of $(-1,1)\cap\Q$ would be a zero-neighborhood on which $g$ is bounded, a contradiction.
\end{proof}

\begin{theorem}[Strict intermediate topology]\label{thm:main}
For every nonzero $f\in D$,
\[
\boxed{\tu\subsetneq\sigma_f\subsetneq\tau_f},
\qquad D_{\sigma_f}=\{0\}.
\]
In particular, $\tu\vee\beta_f\ne\tau_f$.
\end{theorem}
\begin{proof}
Apply Theorem~\ref{thm:vanishing} to $\mathcal S=\{f_a:a\in\R\}$. The function $f$ is not $\sigma_f$-continuous, but is $\tau_f$-continuous. Proposition~\ref{prop:basic} gives the other inclusions and their stated strictness.
\end{proof}

\begin{corollary}[Failure of representation by rational-valued observations]\label{cor:representation}
The topology $\sigma_f$ cannot be generated by $\tu$ and any family of maps in $D$, with the target $\Q$ carrying its usual topology. Nevertheless $\sigma_f$ is a nondiscrete Hausdorff TVG topology strictly stronger than usual and $1/n\to0$ in $\sigma_f$.
\end{corollary}
\begin{proof}
Every member of a proposed generating family would be $\sigma_f$-continuous and hence zero by Theorem~\ref{thm:main}. Such a family generates only $\tu$, contradicting strictness. The remaining assertions follow from Proposition~\ref{prop:basic}.
\end{proof}
This conclusion concerns the specified real-valued (equivalently, usual-rational-valued) coefficient construction. The topology $\sigma_f$ itself is built from a Hamel function, with circle-valued observations. Thus it would be misleading to call it ``not Hamel-based'' without specifying the target and the construction.

\section{Connectedness, characters, and additive lifts}
\begin{proposition}\label{prop:disconnected}
Each of $(\R,\beta_f)$, $(\R,\sigma_f)$, and $(\R,\tau_f)$ is totally disconnected. In each case the continuous characters $p\circ f_a$ separate points.
\end{proposition}
\begin{proof}
A nonzero rational-linear map $f:\R\to\Q$ is onto. If $x\ne y$, choose $z$ with $f(z)=1/2$ and set $a=z/(x-y)$. Then
\[
f_a(x)-f_a(y)=1/2,
\]
so $p(f_a(x))\ne p(f_a(y))$.

The subgroup $\Q/\Z$ of $\T$ is totally disconnected. Indeed, remove a point $p(\theta)$ with $\theta$ irrational. The remaining circle is homeomorphic to an open real interval, and the image of $\Q/\Z$ is countable. A countable subset of the real line contains no nontrivial interval and hence has only singleton connected subsets. Every connected subset of any of our three spaces has singleton image under every $p\circ f_a$. Separation of points therefore forces that subset to be a singleton.
\end{proof}

\begin{lemma}[Uniqueness of additive lifts]\label{lem:lift}
If additive maps $g,h:\R\to\R$ satisfy $p\circ g=p\circ h$, then $g=h$.
\end{lemma}
\begin{proof}
The difference $g-h$ is an additive map into $\Z$ and is zero by the divisibility argument in Lemma~\ref{lem:noncontinuous}.
\end{proof}
This is the standard uniqueness fact in \cite[Lemma~3.1]{DT2008}; it is included to make the following consequence self-contained.

\begin{corollary}[Continuous characters without continuous additive lifts]\label{cor:lifts}
For $a\ne0$, the $\sigma_f$-continuous character $p\circ f_a$ has a unique additive lift to $\R$, namely $f_a$. That lift is not continuous. The same assertion holds for the domain topology $\beta_f$.
\end{corollary}
\begin{proof}
The map $f_a$ is nonzero and is an additive lift. Uniqueness follows from Lemma~\ref{lem:lift}. It is not continuous as a map into $\Q_{\mathrm{usual}}$ by Theorem~\ref{thm:main}. Since $f_a$ takes values in $\Q$, this is equivalent to failure of continuity as a real-valued function. Indeed, every open subset of $\Q_{\mathrm{usual}}$ has the form $U\cap\Q$ for a usual-open set $U\subseteq\R$, and $f_a^{-1}(U\cap\Q)=f_a^{-1}(U)$. Conversely, continuity into $\Q_{\mathrm{usual}}$ implies real-valued continuity by composition with the continuous inclusion $\Q_{\mathrm{usual}}\hookrightarrow\R_{\mathrm{usual}}$. Since $\beta_f\subseteq\sigma_f$, it is not $\beta_f$-continuous either.
\end{proof}
The lifts here are required to be additive, not necessarily real-linear. Dom\'inguez and Tarieladze's character lifting property requires real-linear lifts, a stronger requirement. Their Example~4.3 uses Hamel functions on the discrete real line: those functions are continuous additive lifts but need not be real-linear. In our construction the unique additive lift itself fails to be continuous. Nonliftable characters on nondiscrete TVGs are not new as a general phenomenon; their Example~4.4 already supplies one on the Bohr real line.

A character is an almost periodic function on any topological group on which it is continuous: its translates are scalar multiples of the character, with scalars in a compact circle. Thus the Hamel characters are almost periodic for $\beta_f$, $\sigma_f$, and $\tau_f$, but are not usual-continuous almost periodic functions. The topology on the domain must be specified.

\section{Comparison with arbitrary TVG topologies}
The following elementary criteria locate a TVG topology relative to $\tu$. They are included as comparison tools, not as novelty claims.

\begin{proposition}\label{prop:criteria}
For a TVG topology $\tau$ on $\R$:
\begin{enumerate}
\item $\tau\subseteq\tu$ if and only if every $\tau$-neighborhood of zero contains an ordinary open interval about zero.
\item $\tu\subseteq\tau$ if and only if some $\tau$-neighborhood of zero is bounded in ordinary absolute value.
\end{enumerate}
Consequently, every zero-neighborhood in a strictly weaker TVG topology is ordinarily unbounded. Every topology weaker than usual is path connected.
\end{proposition}
\begin{proof}
The first equivalence is continuity at zero of the identity from $\R_{\mathrm{usual}}$ into $(\R,\tau)$; for group homomorphisms continuity at zero implies continuity everywhere. For the forward direction in the second equivalence, the interval $(-1,1)$ is a $\tau$-neighborhood. Conversely, let $U\subseteq[-M,M]$ be a $\tau$-neighborhood, choosing $M>0$. For $\varepsilon>0$ choose $c>0$ with $cM<\varepsilon$. The set $cU$ is a $\tau$-neighborhood contained in $(-\varepsilon,\varepsilon)$. This proves continuity of the identity into $\R_{\mathrm{usual}}$ and hence the second equivalence.

A bounded neighborhood in a strictly weaker topology would imply the reverse inclusion and equality, a contradiction. If $\tau\subseteq\tu$, the usual path $s\mapsto(1-s)x+sy$ from $[0,1]$ remains continuous into $(\R,\tau)$ for every pair $x,y$.
\end{proof}
The boundedness in this proposition is numerical boundedness, not an unspecified notion of boundedness for TVGs. Combining its two criteria characterizes $\tu$ itself.

For every strictly stronger TVG topology and $v\ne0$, the orbit map $a\mapsto av$ from the usual scalar line is discontinuous at zero. Otherwise composition with the fixed homeomorphism $x\mapsto x/v$ would make the identity from $\R_{\mathrm{usual}}$ to the stronger topology continuous, forcing equality. Requiring joint scalar continuity instead gives an ordinary Hausdorff real topological vector space, whose one-dimensional topology must be the usual one.

\section{Families, common refinements, and the resulting diagram}
A family $\mathcal S\subseteq D$ is \emph{invariant under rescaling} if $s_a\in\mathcal S$ for every $s\in\mathcal S$ and $a\in\R$. For such a family, let $\tau_{\mathcal S}$ be generated by $\tu$ and the maps $s\in\mathcal S$ into $\Q_{\mathrm{usual}}$.

\begin{proposition}\label{prop:family}
Suppose that $\mathcal S\subseteq D$ is invariant under rescaling and contains a nonzero function. Then $\sigma_{\mathcal S}$ and $\tau_{\mathcal S}$ are Hausdorff nondiscrete TVG topologies, both totally disconnected, and
\[
\tu\subsetneq\sigma_{\mathcal S}\subsetneq\tau_{\mathcal S},
\qquad D_{\sigma_{\mathcal S}}=\{0\}.
\]
In particular, taking $\mathcal S=D$ gives common nondiscrete refinements of all the respective single-function constructions:
\[
\bigvee_{0\ne f\in D}\sigma_f=\sigma_D,
\qquad
\bigvee_{0\ne f\in D}\tau_f=\tau_D,
\qquad \sigma_D\subsetneq\tau_D.
\]
\end{proposition}
\begin{proof}
Rescaling invariance proves continuity of fixed dilations, and the identity coordinate proves Hausdorffness. Every basic neighborhood uses finitely many observations; $s(1/n)=s(1)/n\to0$ for each, so $1/n\to0$ in both topologies. A nonzero $s\in\mathcal S$ supplies a usual-discontinuous character, proving the first strict inclusion. Theorem~\ref{thm:vanishing} makes this same $s$ discontinuous in $\sigma_{\mathcal S}$, whereas it is continuous in $\tau_{\mathcal S}$, proving the second. The characters $p\circ s_a$ separate points exactly as in Proposition~\ref{prop:disconnected}, proving total disconnectedness. For $\mathcal S=D$, every defining function or character occurs in a single-function construction; conversely all such constructions use members of $D$. This proves both supremum identities.
\end{proof}
For a fixed nonzero $f$, our conclusions give the two strict chains
\[
\begin{aligned}
t&\subsetneq\tu\subsetneq\sigma_f\subsetneq\tau_f,\\
t&\subsetneq\beta_f\subsetneq\sigma_f\subsetneq\tau_f,
\end{aligned}
\]
with $\beta_f$ and $\tu$ incomparable. Strictness of $\beta_f\subsetneq\sigma_f$ follows because $\sigma_f$ contains $\tu$, which is not contained in $\beta_f$. This is a branching order of topologies rather than a single universal chain. The earlier topologies $\rho_{(r_n)}$ of \cite{Y2026} can be inserted between $t$ and $\tu$ in the first chain.

For comparison, the companion manuscript \cite{YHamel} proves
\[
D_{\tau_f}=\{f_a:a\in\R\},\qquad
\tau_g\subseteq\tau_f\ \Longleftrightarrow\ g=f_b\text{ for some }b\ne0
\ \Longleftrightarrow\ \tau_g=\tau_f.
\]
These are recalled results from that manuscript; their proofs are not required for any theorem above. They show why the vanishing theorem for $\sigma_f$ is a substantial change in the behavior of rational-valued observations, despite $\sigma_f$ lying between $\tu$ and $\tau_f$.

\section{Completion models}\label{sec:completion}
We identify the group completions using standard Pontryagin duality for locally compact abelian groups. These arguments do not assume that our original groups are locally compact. For an abelian topological group $G$, its Pontryagin dual $\widehat G$ is the group of continuous characters $G\to\T$, with the compact-open topology when a topology on the dual is needed. A completion below is the Hausdorff group completion for the two-sided uniformity, not the Stone--Cech compactification.

We use the following standard facts. The dual of a discrete abelian group is compact; evaluation identifies an LCA group with its bidual. A subgroup of an LCA group is dense precisely when its annihilator in the dual is zero. A compact abelian group is connected precisely when its discrete dual is torsion-free, and it is metrizable precisely when that dual is countable. These facts are part of classical Pontryagin duality; see \cite[Chapters~5--7]{Morris1977}. The subgroup/quotient duality is also summarized in \cite{Elkies2002}; the annihilator criterion and the connectedness criterion are recorded in \cite[Appendix~C, Theorem~C.13 and the paragraph following it]{EW2010}. Continuous homomorphisms into a complete Hausdorff group extend uniquely from a dense subgroup to its group completion.

Let $\R_d$ denote the additive real line with the discrete topology, and set
\[
K=\widehat{\R_d}=\operatorname{Hom}(\R_d,\T).
\]
All algebraic characters of $\R_d$ are continuous. We give $K$ the topology of pointwise convergence, which here is its compact-open topology. By duality, every continuous character of $K$ has the form
\[
\operatorname{ev}_a(\kappa)=\kappa(a)\qquad(a\in\R),
\]
with unique parameter $a$. The usual Bohr compactification of the real line is this $K$, with dense evaluation map $e$ defined below. In particular, $K$ is fixed independently of $f$.

For nonzero $f\in D$ define homomorphisms $e,j_f:\R\to K$ by
\[
e(x)(b)=p(bx),\qquad
j_f(x)(a)=p(f(ax)).
\]
For each $x$, these are characters of the parameter group $\R_d$. Define
\[
I_f(x)=(x,j_f(x))\in\R_{\mathrm{usual}}\times K,\qquad
J_f(x)=(e(x),j_f(x))\in K\times K.
\]
The first coordinate makes $I_f$ injective, and ordinary characters make $J_f$ injective. By their defining initial topologies, $I_f$ is a topological group embedding of $(\R,\sigma_f)$, and $J_f$ is a topological group embedding of $(\R,\beta_f)$.

\begin{lemma}[No mixed additive relation]\label{lem:mixed}
If $c,a\in\R$ satisfy
\[
p(cx+f(ax))=0\qquad\text{for every }x\in\R,
\]
then $c=a=0$.
\end{lemma}
\begin{proof}
The additive function $x\mapsto cx+f(ax)$ takes integer values, so it is zero by Lemma~\ref{lem:lift} and its divisibility proof. It follows that $cx=-f(ax)\in\Q$ for every real $x$. If $c\ne0$, then $c\R=\R$ cannot be contained in $\Q$. Hence $c=0$. If $a\ne0$, multiplication by $a$ is onto, and $f(ax)=0$ for every $x$ would force $f=0$. Thus $a=0$.
\end{proof}

\begin{theorem}[Explicit completions]\label{thm:completions}
The images of $I_f$ and $J_f$ are dense in their respective target groups. Consequently
\[
\boxed{
\operatorname{comp}(\R,\sigma_f)\cong\R_{\mathrm{usual}}\times K,\qquad
\operatorname{comp}(\R,\beta_f)\cong K\times K.
}
\]
These are isomorphisms of topological groups, with the specified embeddings $I_f$ and $J_f$.
\end{theorem}
\begin{proof}
Every continuous character of $\R_{\mathrm{usual}}\times K$ has the form
\[
(r,\kappa)\longmapsto p(cr)+\kappa(a)\qquad(c,a\in\R).
\]
Indeed, a homomorphism on a product is the sum of its restrictions to the factors, and the characters of the two factors are as described above. Its restriction to $I_f(\R)$ is $x\mapsto p(cx+f(ax))$. Lemma~\ref{lem:mixed} says that the only character annihilating $I_f(\R)$ is the zero character. Since the target is LCA, the annihilator criterion proves density.

Similarly, every character of $K\times K$ has the form
\[
(\kappa_1,\kappa_2)\longmapsto\kappa_1(c)+\kappa_2(a).
\]
Its restriction to $J_f(\R)$ is again $p(cx+f(ax))$, so the same criterion proves density. The target $K\times K$ is compact and complete; $\R_{\mathrm{usual}}\times K$ is also complete. A dense topological group embedding into a complete Hausdorff group identifies that target with the group completion.
\end{proof}

\begin{remark}
The completion models do not depend on the choice of nonzero $f$, but the embeddings do. Thus isomorphic completions do not imply equality of the topologies on the fixed real line. Each fixed scalar dilation also extends continuously to these completions. If
\[
(L_c\kappa)(a)=\kappa(ca),
\]
the extensions are $(r,\kappa)\mapsto(cr,L_c\kappa)$ for $I_f$ and $(\kappa_1,\kappa_2)\mapsto(L_c\kappa_1,L_c\kappa_2)$ for $J_f$. These actions satisfy the real vector space identities. No joint continuity with the usual scalar topology is asserted.
\end{remark}

\section{All continuous characters and real-valued additive functions}\label{sec:all-duals}
\begin{theorem}[Complete description of continuous characters]\label{thm:characters}
For either $\eta_f=\sigma_f$ or $\eta_f=\beta_f$, the continuous characters are exactly
\[
\boxed{
\widehat{(\R,\eta_f)}
=
\{\chi_{c,a}:c,a\in\R\},\qquad
\chi_{c,a}(x)=p(cx+f(ax)).
}
\]
The parameters $c,a$ are unique. This equality describes characters as functions, and hence as an abstract group; it makes no identification of the compact-open topologies of the two duals.
\end{theorem}
\begin{proof}
Every displayed character is continuous for $\beta_f$, and hence for $\sigma_f$, by construction. Conversely, a continuous character of either original group extends uniquely to its completion, since $\T$ is complete and homomorphisms are uniformly continuous. The descriptions of characters of the two completion models in the proof of Theorem~\ref{thm:completions} yield the formula. If $\chi_{c,a}=\chi_{d,b}$, then
\[
p((c-d)x+f((a-b)x))=0
\]
for every $x$. Lemma~\ref{lem:mixed} gives $c=d$ and $a=b$.
\end{proof}

\begin{corollary}[The Bohr modification]\label{cor:bohr}
The topologies $\beta_f$ and $\sigma_f$ have the same continuous characters. Moreover, $\beta_f$ is the Bohr modification of $\sigma_f$, meaning the coarsest topology making all its continuous characters continuous.
\end{corollary}
\begin{proof}
The character sets agree by Theorem~\ref{thm:characters}. Every $\chi_{c,a}$ is a sum of two defining characters of $\beta_f$, and those defining characters occur among the $\chi_{c,a}$. Their initial topology is therefore exactly $\beta_f$.
\end{proof}
In the terminology of group duality, these are distinct compatible topologies: they have the same continuous character group. This is compatible with the strict inclusion $\beta_f\subsetneq\sigma_f$ proved earlier.

\begin{theorem}[All continuous real-valued additive functions]\label{thm:real-duals}
With the real target carrying its usual topology,
\[
\boxed{
\operatorname{Hom}_{\mathrm{cont}}((\R,\sigma_f),\R_{\mathrm{usual}})
=\{x\mapsto cx:c\in\R\},
}
\]
whereas
\[
\boxed{
\operatorname{Hom}_{\mathrm{cont}}((\R,\beta_f),\R_{\mathrm{usual}})
=\{0\}.
}
\]
\end{theorem}
\begin{proof}
A continuous additive function extends uniquely to a continuous homomorphism from the appropriate completion into the complete group $\R_{\mathrm{usual}}$. A continuous homomorphism from a compact group into the additive real line is zero: its image is a compact subgroup, and any subgroup containing a nonzero real number contains its unbounded integer multiples.

For $\beta_f$, the completion $K\times K$ is compact, so the extension is zero. For $\sigma_f$, an extension on $\R_{\mathrm{usual}}\times K$ vanishes on $\{0\}\times K$. Its restriction to the real factor is a usual-continuous additive function, hence $r\mapsto cr$. Restricting along $I_f(x)=(x,j_f(x))$ gives $x\mapsto cx$. Conversely these scalar multiples of the identity are continuous in $\sigma_f$.
\end{proof}
In particular, the full real-valued additive dual of $\sigma_f$ is exactly the usual one, although $\sigma_f$ is strictly finer and has additional continuous characters. This strengthens the earlier vanishing statement for rational-valued additive functions. No assertion about all continuous nonadditive real functions is intended.

\begin{corollary}[Exactly which characters lift continuously]\label{cor:all-lifts}
Every $\chi_{c,a}$ has a unique additive real-valued lift
\[
h_{c,a}(x)=cx+f(ax).
\]
For the topology $\sigma_f$, this lift is continuous if and only if $a=0$. For the topology $\beta_f$, it is continuous if and only if $c=a=0$.
\end{corollary}
\begin{proof}
The displayed map is an additive lift, and Lemma~\ref{lem:lift} proves uniqueness. In $\sigma_f$, Theorem~\ref{thm:real-duals} says that a continuous lift must be $dx$. Then $(c-d)x+f(ax)=0$, which by the argument of Lemma~\ref{lem:mixed} forces $a=0$ and $c=d$. The converse is immediate. In $\beta_f$ every continuous real-valued additive function is zero, so only the trivial character has a continuous additive lift; uniqueness of parameters gives $c=a=0$.
\end{proof}
Thus we distinguish three notions explicitly: an algebraic additive lift, a continuous additive lift, and a continuous real-linear lift. The first exists for every continuous character in these examples; the second has been completely classified above.

\section{Equality, inclusion, and multiplicity}\label{sec:equality}
\begin{theorem}[Exact comparison of the modular topologies]\label{thm:modular-comparison}
Let $f,g\in D\setminus\{0\}$. The following are equivalent:
\begin{enumerate}
\item $\sigma_g\subseteq\sigma_f$;
\item $\sigma_g=\sigma_f$;
\item $g(x)=f(bx)$ for all $x$, for some $b\in\R\setminus\{0\}$;
\item $\beta_g=\beta_f$;
\item $\beta_g\subseteq\beta_f$.
\end{enumerate}
Consequently distinct members within either modular family are incomparable.
\end{theorem}
\begin{proof}
Suppose $\sigma_g\subseteq\sigma_f$. The character $p\circ g$ is $\sigma_g$-continuous and therefore $\sigma_f$-continuous. By Theorem~\ref{thm:characters}, there exist $c,b\in\R$ such that
\[
p(g(x))=p(cx+f(bx))\qquad(x\in\R).
\]
Uniqueness of additive lifts gives $g(x)=cx+f(bx)$. Since both $g$ and $f_b$ are rational-valued, $cx$ is rational for every real $x$, forcing $c=0$. Also $b\ne0$, since $g\ne0$. This proves (1)$\Rightarrow$(3).

If $g=f_b$ with $b\ne0$, then $g_a=f_{ba}$. As $a$ ranges over $\R$, so does $ba$. Hence the defining character families coincide for both $\beta$ and $\sigma$, proving (3)$\Rightarrow$(2),(4). Equality implies the respective inclusion. Finally, starting with (5), the same argument using the character classification for $\beta_f$ gives (3). All conditions are equivalent. Either inclusion between two distinct members would imply equality, proving incomparability.
\end{proof}

\begin{corollary}[Number of distinct modular topologies]\label{cor:modular-count}
Each family $\{\sigma_f:0\ne f\in D\}$ and $\{\beta_f:0\ne f\in D\}$ contains exactly $2^{\mathfrak c}$ distinct topologies on the fixed real line, where $\mathfrak c=|\R|$.
\end{corollary}
\begin{proof}
The rational dimension of $\R$ is $\mathfrak c$. A rational-linear functional is determined by arbitrary rational values on a Hamel basis, so
\[
|D|=|\Q|^{\mathfrak c}=2^{\mathfrak c}.
\]
The lower bound uses only values zero and one; the upper bound follows from $|\Q|\le2^{\aleph_0}$. For fixed nonzero $f$, the functions $f_b$, $b\ne0$, are pairwise distinct: equality for $b\ne d$ would make $f((b-d)x)=0$ for every $x$, contrary to $f\ne0$. By Theorem~\ref{thm:modular-comparison}, each equality class therefore has cardinality $\mathfrak c$. Partitioning a set of cardinality $2^{\mathfrak c}$ into classes of size $\mathfrak c$ gives $2^{\mathfrak c}$ classes: if their number is $\lambda$, cardinal multiplication gives $2^{\mathfrak c}=\lambda\mathfrak c=\max\{\lambda,\mathfrak c\}$, and Cantor's inequality forces $\lambda=2^{\mathfrak c}$.
\end{proof}
This counts distinct collections of open sets on the same underlying real line. It does not classify homeomorphism types. In fact the completion models in Theorem~\ref{thm:completions} are independent of $f$.

\section{First countability, completeness, and local compactness}\label{sec:countability}
\begin{theorem}\label{thm:nonmetric}
Neither $\sigma_f$ nor $\beta_f$ is first countable or metrizable. Neither group is complete or locally compact. Both completions are connected and nonmetrizable; the completion of $\beta_f$ is compact, whereas the completion of $\sigma_f$ is locally compact and noncompact.
\end{theorem}
\begin{proof}
The discrete dual of $K$ is $\R_d$, which is uncountable and torsion-free. Standard compact/discrete duality therefore says that $K$ is nonmetrizable and connected. In particular it is nontrivial. The completion models $\R_{\mathrm{usual}}\times K$ and $K\times K$ are connected and contain a topological copy of $K$, so neither is metrizable. The second is compact. The first is locally compact, since $\R_{\mathrm{usual}}$ is locally compact and $K$ compact, and is noncompact because it maps continuously onto the usual real line.

If either original group were first countable, the Birkhoff--Kakutani theorem would give a compatible invariant metric; invariance can be taken two-sided here because the groups are abelian. Its metric completion realizes the group completion and is metrizable, contradicting the preceding paragraph. A metrizable space is first countable, so both possibilities are excluded. For the metrizability theorem and the compact-dual countability criterion, see \cite[Theorem~8.3 and Remark~6.11]{MegNotes}.

Proposition~\ref{prop:disconnected} makes both original groups totally disconnected, while their completions are nontrivial and connected. Thus their dense images are proper. A complete subgroup of a Hausdorff complete group is closed, so a proper dense image cannot be complete. A locally compact subgroup of a Hausdorff topological group is also closed; see \cite[Theorem~3.7]{MegNotes}. Hence these proper dense subgroups are not locally compact.
\end{proof}
This supplies explicit totally disconnected dense subgroups of connected completions. Total disconnectedness need not be preserved by group completion. Nor does a locally compact completion make its dense subgroup locally compact.

\section{Continuous homomorphisms and group isomorphisms}\label{sec:rigidity}
The comparison theorem concerns the identity map on the fixed real line. We now allow every additive map between two members of the modular families.

\begin{theorem}[Rigidity of continuous homomorphisms]\label{thm:rigidity}
Let $f,g\in D\setminus\{0\}$. Choose independently $\eta_f\in\{\beta_f,\sigma_f\}$ and $\eta_g\in\{\beta_g,\sigma_g\}$. Every continuous additive map
\[
T:(\R,\eta_f)\longrightarrow(\R,\eta_g)
\]
has the form $T(x)=cx$ for a real constant $c$. If $c\ne0$, continuity holds exactly when $\eta_g\subseteq\eta_f$.
\end{theorem}
\begin{proof}
Every ordinary character is continuous on the target. Consequently $p\circ T$ is continuous on the source. By Theorem~\ref{thm:characters} and uniqueness of additive lifts,
\[
T(x)=cx+f(ax)
\]
for some $c,a\in\R$. Fix an irrational number $r$. Applying the same argument to $p(rT(x))$ gives $rT(x)=dx+f(bx)$ for suitable $d,b$. Hence
\[
(rc-d)x=f(bx)-rf(ax).
\]
The right side takes values in the countable set $\Q+r\Q$. A nonzero real multiple of the identity has uncountable image, so $rc=d$. We obtain $rf_a=f_b$. If $a\ne0$, then $f_a$ is onto $\Q$; evaluating at a point where $f_a=1$ contradicts the rational-valuedness of $f_b$. Thus $a=0$ and $T=c\,\mathrm{id}$.

For $c\ne0$, multiplication by $c$ is a homeomorphism of the source TVG. Composing its inverse with $T$ shows that continuity of $T$ is equivalent to continuity of the identity from $\eta_f$ to $\eta_g$, which is the stated inclusion.
\end{proof}

\begin{corollary}[Isomorphism classification]\label{cor:isomorphism}
Within either the $\beta$-family or the $\sigma$-family, two groups are topologically isomorphic if and only if their topologies are equal, equivalently $g=f_b$ for some $b\ne0$. Between distinct members of either family, the only continuous homomorphism is zero. Every continuous endomorphism of a member is a scalar multiplication. Each family therefore contains exactly $2^{\mathfrak c}$ pairwise nonisomorphic topological groups.
\end{corollary}
\begin{proof}
A nonzero continuous homomorphism must be a nonzero scalar multiplication by Theorem~\ref{thm:rigidity}. It forces inclusion, and Theorem~\ref{thm:modular-comparison} turns inclusion within either family into equality. Conversely, equal topologies admit every scalar multiplication, with every nonzero one a homeomorphism. Apply Corollary~\ref{cor:modular-count} for the number of classes.
\end{proof}
This classifies topological group isomorphisms, not arbitrary homeomorphisms. In particular, each family gives pairwise nonisomorphic dense groups with the same completion model.

\begin{corollary}[Maps between the two families]
There is no nonzero continuous homomorphism from $(\R,\beta_f)$ to $(\R,\sigma_g)$. A nonzero continuous homomorphism from $(\R,\sigma_f)$ to $(\R,\beta_g)$ exists exactly when $g=f_b$ for some $b\ne0$; in that case all scalar multiplications are continuous.
\end{corollary}
\begin{proof}
In the first direction Theorem~\ref{thm:rigidity} would imply $\tu\subseteq\sigma_g\subseteq\beta_f$, contradicting Proposition~\ref{prop:basic}. In the second direction, inclusion $\beta_g\subseteq\sigma_f$ makes $p\circ g$ continuous on $\sigma_f$. The proof of Theorem~\ref{thm:modular-comparison} gives $g=f_b$. Conversely this relation implies $\beta_g=\beta_f\subseteq\sigma_f$.
\end{proof}

\section{Compact sets and convergent sequences}\label{sec:compact}
We use Glicksberg's theorem: a locally compact abelian group and its Bohr topology have exactly the same compact subsets \cite[Theorem~1.2]{Glicksberg1962}. Its application here concerns the locally compact completion model, not an assumption that our dense groups are locally compact.

\begin{theorem}[Compact subsets and sequences of the modular pair]\label{thm:compact}
The topologies $\beta_f$ and $\sigma_f$ have the same compact subsets and induce the same topology on each of them. They have the same convergent sequences with the same limits. Their common character group also has the same compact-open topology for both domain topologies.
\end{theorem}
\begin{proof}
Put $L=\R_{\mathrm{usual}}\times K$, and give it also its Bohr topology $L^+$, generated by its continuous characters. Their restrictions under $I_f$ are exactly the functions $\chi_{c,a}$ from Theorem~\ref{thm:characters}. Thus $I_f$ embeds $(\R,\beta_f)$ into $L^+$ and $(\R,\sigma_f)$ into $L$.

If $C$ is $\beta_f$-compact, then $I_f(C)$ is compact in $L^+$. Glicksberg's theorem makes it compact in $L$, hence $C$ is $\sigma_f$-compact. The converse follows from $\beta_f\subseteq\sigma_f$. On a common compact subset, the identity from the finer topology to the coarser Hausdorff topology is a homeomorphism.

A convergent sequence together with its limit is compact. Equality of the induced topologies on compact subsets therefore gives equality of convergent sequences and limits in both directions. Finally, compact-open neighborhoods in the common character set use precisely these compact subsets and their character values. These data agree, proving the last assertion.
\end{proof}

\begin{corollary}[Explicit convergence tests]\label{cor:sequences}
For a sequence $(x_n)$ and $x\in\R$, convergence in either $\beta_f$ or $\sigma_f$ is equivalent to the two conditions
\[
x_n\longrightarrow x\text{ usually},\qquad
\|f(a(x_n-x))\|_{\T}\longrightarrow0\quad(a\in\R).
\]
Convergence in $\tau_f$ is equivalent to usual convergence and
\[
f(a(x_n-x))\longrightarrow0\text{ in }\Q_{\mathrm{usual}}
\quad(a\in\R).
\]
\end{corollary}
\begin{proof}
The initial-topology definitions give the tests for $\sigma_f$ and $\tau_f$. Theorem~\ref{thm:compact} transfers the $\sigma_f$ test to $\beta_f$.
\end{proof}
Thus sequences cannot distinguish $\beta_f$ from its strictly finer topology $\sigma_f$. This does not imply that all nets have the same convergence.

\begin{proposition}[Compact-set tests]\label{prop:compact-tests}
A subset $C\subseteq\R$ is compact for $\beta_f$ or $\sigma_f$ if and only if it is usual-compact and each restriction $p\circ f_a|_C$ is usual-continuous. It is $\tau_f$-compact if and only if it is usual-compact and each restriction $f_a|_C:C\to\Q_{\mathrm{usual}}$ is usual-continuous.
\end{proposition}
\begin{proof}
For $\sigma_f$-compact $C$, its continuous identity map onto the usual Hausdorff subspace $C$ is a homeomorphism. The restricted observations are then usual-continuous. Conversely, if the stated usual-continuity conditions hold, the map from usual $C$ into $(\R,\sigma_f)$ is continuous by the initial-topology definition; its image is compact. Theorem~\ref{thm:compact} handles $\beta_f$. The same argument with real-valued observations proves the assertion for $\tau_f$.
\end{proof}

\begin{theorem}[A universal separating sequence]\label{thm:universal-sequence}
The sequence
\[
x_n=n!(3-2\sqrt2)^{n^2}\qquad(n\ge1)
\]
converges to zero in $\sigma_f$ and $\beta_f$ for every nonzero $f\in D$, but does not converge to zero in $\tau_f$ for any such $f$. More precisely, for every fixed $a$ the value $f(ax_n)$ is eventually an integer, while for a suitable $a$ its absolute value tends to infinity.
\end{theorem}
\begin{proof}
Define positive integers $P_n,Q_n$ by
\[
(3+2\sqrt2)^{n^2}=P_n+Q_n\sqrt2.
\]
The binomial expansion gives integer coefficients, and conjugation gives
\[
x_n=n!(P_n-Q_n\sqrt2).
\]
Since $0<3-2\sqrt2<1$, exponential decay in $n^2$ dominates $n!$, so $x_n\to0$ usually. For every fixed $a$,
\[
f(ax_n)=n!\bigl(P_nf(a)-Q_nf(a\sqrt2)\bigr).
\]
Both coefficient values are rational. Their fixed denominators divide $n!$ eventually, so $f(ax_n)\in\Z$ eventually. The modular convergence test in Corollary~\ref{cor:sequences} proves convergence in $\sigma_f$ and $\beta_f$.

Choose $a$ with $f(a)=1$, which is possible because $f$ is onto $\Q$, and put $r=f(a\sqrt2)\in\Q$. Then
\[
f(ax_n)=n!Q_n(\sqrt2-r)+x_n.
\]
Here $\sqrt2-r\ne0$ and $n!Q_n\to\infty$. Hence $|f(ax_n)|\to\infty$, excluding $\tau_f$-convergence.
\end{proof}
The same numerical sequence works simultaneously for every nonzero rational-valued coefficient functional. It also converges in the common modular refinement $\sigma_D$ of Section~8.

\begin{corollary}[Different compact subsets]\label{cor:compact-separation}
The set $C=\{0,x_1,x_2,\ldots\}$ from Theorem~\ref{thm:universal-sequence} is compact in $\beta_f$ and $\sigma_f$ but not in $\tau_f$.
\end{corollary}
\begin{proof}
A convergent sequence together with its limit is compact. If $C$ were $\tau_f$-compact, its continuous identity map onto the usual Hausdorff subspace $C$ would be a homeomorphism. Usual convergence $x_n\to0$ would then imply $\tau_f$-convergence, a contradiction.
\end{proof}

\begin{theorem}[Sequences in a finite rational span]\label{thm:finite-span-sequences}
Let $b_1,\ldots,b_m\in\R$ be rationally independent and write
\[
x_n=\sum_{j=1}^m q_{n,j}b_j,\qquad q_{n,j}\in\Q.
\]
For every nonzero $f\in D$, convergence $x_n\to0$ in $\tau_f$ is equivalent to $q_{n,j}\to0$ usually for each $j$. Convergence in $\sigma_f$, equivalently in $\beta_f$, is equivalent to usual convergence of $x_n$ together with
\[
\left\|\frac{q_{n,j}}d\right\|_{\T}\longrightarrow0
\qquad(1\le j\le m,\ d\in\mathbb N,\ d\ge1).
\]
Thus all three convergence tests on a fixed finite rational span are independent of the choice of nonzero $f$.
\end{theorem}
\begin{proof}
The rational-linear map
\[
a\longmapsto(f(ab_1),\ldots,f(ab_m))
\]
is onto $\Q^m$. Otherwise its image would have a nonzero rational annihilating vector $(r_j)$, giving
\[
f\left(a\sum_{j=1}^m r_jb_j\right)=0\qquad(a\in\R).
\]
Independence makes the scalar in parentheses nonzero, contradicting $f\ne0$.

If $x_n\to0$ in $\tau_f$, choose an $a$ whose displayed image is the $j$th coordinate vector. Then $f(ax_n)=q_{n,j}\to0$. Conversely, convergence of every coordinate implies usual convergence of $x_n$ and real convergence of each finite sum $f(ax_n)=\sum_jq_{n,j}f(ab_j)$.

For $\sigma_f$-convergence, choose the image vector with $j$th entry $\frac1d$ and the others zero. This gives the necessary modular conditions. Conversely, each fixed rational coefficient $f(ab_j)$ is an integer divided by some positive integer $d$. The stated conditions and continuity of addition in $\T$ give modular convergence of $\sum_jq_{n,j}f(ab_j)$. Apply Corollary~\ref{cor:sequences}.
\end{proof}

\begin{corollary}[Least rational dimension for separation]
A sequence contained in a one-dimensional rational subspace cannot converge to zero in $\sigma_f$ without converging to zero in $\tau_f$. Dimension two suffices.
\end{corollary}
\begin{proof}
In dimension one, write $x_n=q_nb$ with $b\ne0$. Usual convergence forces $q_n\to0$, so Theorem~\ref{thm:finite-span-sequences} gives $\tau_f$-convergence. The universal sequence belongs to $\operatorname{span}_{\Q}\{1,\sqrt2\}$, proving sufficiency of dimension two.
\end{proof}
Theorem~\ref{thm:finite-span-sequences} classifies sequences in finite rational spans. For sequences with infinite-dimensional rational span, Corollary~\ref{cor:sequences} remains the general test; no further classification is asserted here.

\section{Character subspaces and all precompact TVG topologies}\label{sec:precompact}
Let
\[
\mathcal A=\operatorname{Hom}((\R,+),\T)
\]
be the group of all algebraic characters, with no continuity assumption. It has a real vector space structure given by pointwise addition and the scalar action
\[
(a\cdot\chi)(x)=\chi(ax).
\]
This uses rescaling of the argument. Multiplication of an arbitrary circle value by an arbitrary real scalar is not being defined. Additivity of $\chi$ verifies $(a+b)\cdot\chi=a\cdot\chi+b\cdot\chi$, and the other vector space identities follow directly.

For a subgroup $H\subseteq\mathcal A$, write $\theta_H$ for its initial character topology. We recall the standard precompact character representation, a consequence of compact Pontryagin duality \cite{Morris1977}: a point-separating character subgroup generates a Hausdorff precompact topology whose full continuous character group is that subgroup. Indeed, embed the group into $\T^H$ and take its compact closure. A continuous character extends to that closure and then to the product; characters of the product are finite integer combinations of coordinate characters. The same representation shows that every Hausdorff precompact abelian group has its original topology generated by its continuous characters.

\begin{theorem}[Classification of precompact TVG topologies]\label{thm:precompact}
The assignment $H\mapsto\theta_H$ is an order-preserving bijection from the nonzero real vector subspaces of $\mathcal A$ onto the Hausdorff precompact TVG topologies on $\R$. In particular,
\[
\theta_H\subseteq\theta_J\quad\Longleftrightarrow\quad H\subseteq J.
\]
\end{theorem}
\begin{proof}
A nonzero scalar-invariant subgroup $H$ separates points. Given $x\ne0$, choose $\chi\in H$ and $y$ with $\chi(y)\ne0$. Then $(y/x)\cdot\chi$ belongs to $H$ and is nonzero at $x$. Its initial topology is therefore Hausdorff and precompact. Closure of $H$ under the scalar action makes every fixed dilation continuous, so this topology is a TVG topology.

Conversely let $\theta$ be a Hausdorff precompact TVG topology, and let $H$ be its continuous character group. The representation recalled above gives $\theta=\theta_H$ and ensures that $H$ separates points and is nonzero. Composition with fixed dilations shows scalar invariance, so $H$ is a real subspace. Finally, inclusion of character subgroups gives inclusion of their initial topologies. The converse follows because every character continuous for a coarser topology is continuous for the finer one, and their full character groups are $H$ and $J$.
\end{proof}

\begin{corollary}[Minimal precompact topologies and absence of a least topology]\label{cor:minimal}
A Hausdorff precompact TVG topology is minimal among Hausdorff TVG topologies if and only if its character subspace has real dimension one. There is no least Hausdorff TVG topology on $\R$.
\end{corollary}
\begin{proof}
Every coarser group topology of a precompact group is precompact. Thus Theorem~\ref{thm:precompact} reduces minimality to having no proper nonzero real subspace, which is equivalent to dimension one.

For the second assertion put
\[
H_0=\{x\mapsto p(cx):c\in\R\},\qquad
V_f=\{p\circ f_a:a\in\R\},\qquad \alpha_f=\theta_{V_f}.
\]
Both $H_0$ and $V_f$ are one-dimensional real subspaces of $\mathcal A$. Injectivity of their scalar parametrizations follows respectively from the usual real characters and Lemma~\ref{lem:lift}. Lemma~\ref{lem:mixed} gives $H_0\cap V_f=\{0\}$. The topology $\theta_{H_0}$ is $t$. A Hausdorff TVG topology coarser than both $t$ and $\alpha_f$ would be precompact, with a nonzero character subspace contained in their zero intersection, a contradiction. Hence these two topologies have no common coarser Hausdorff TVG topology; a least topology cannot exist.
\end{proof}
This preserves the distinction between minimality and being least. In particular, the author's minimality result for $t$ \cite{Y1981} does not imply that $t$ is contained in every Hausdorff TVG topology. The supplementary topology $\alpha_f$ also gives the identities
\[
\beta_f=t\vee\alpha_f,\qquad \sigma_f=\tu\vee\alpha_f.
\]
The corollary classifies minimal topologies within the precompact class. It does not assert that every minimal Hausdorff TVG topology is precompact.

\begin{theorem}[Number of precompact topologies]\label{thm:precompact-count}
There are exactly $2^{2^{\mathfrak c}}$ Hausdorff precompact TVG topologies on the fixed real line. Exactly $2^{\mathfrak c}$ of these are minimal Hausdorff TVG topologies.
\end{theorem}
\begin{proof}
The upper bound $|\mathcal A|\le|\T|^{|\R|}=2^{\mathfrak c}$ is immediate. The map $D\to\mathcal A$, $f\mapsto p\circ f$, is injective by Lemma~\ref{lem:lift}, so $|\mathcal A|=2^{\mathfrak c}$. An infinite-dimensional vector space over a field of cardinality $\mathfrak c$ has cardinality equal to the maximum of $\mathfrak c$ and its basis cardinality. Since $2^{\mathfrak c}>\mathfrak c$, the real dimension of $\mathcal A$ is $2^{\mathfrak c}$.

Spans of distinct subsets of a basis are distinct subspaces, giving $2^{2^{\mathfrak c}}$ subspaces, and this is also the upper bound supplied by the power set of $\mathcal A$. Omitting the zero subspace does not change the cardinal. Theorem~\ref{thm:precompact} proves the first count. Each real line has $\mathfrak c$ nonzero vectors, so the $2^{\mathfrak c}$ nonzero vectors of $\mathcal A$ partition into $2^{\mathfrak c}$ one-dimensional subspaces. Apply Corollary~\ref{cor:minimal}.
\end{proof}
These counts concern topologies on the fixed set. Unlike Corollary~\ref{cor:isomorphism}, this theorem does not count isomorphism types.

\section{Compatible topologies and the Mackey question}\label{sec:compatible}
Put
\[
H_f=H_0+V_f=\{\chi_{c,a}:c,a\in\R\}.
\]
Two Hausdorff abelian group topologies are called compatible if they have the same continuous characters as functions on the underlying group. This definition does not require any preliminary identification of their compact-open dual topologies.

For precision, a subset $U$ of an abelian topological group is quasi-convex if every $x\notin U$ can be separated from $U$ by a continuous character: that character sends $U$ into the arc $[-\frac14,\frac14]+\Z$ and sends $x$ outside that arc. A topology is locally quasi-convex if it has a zero-neighborhood base of such sets. Let $\mathcal C_f$ denote the compatible locally quasi-convex Hausdorff TVG topologies having character group $H_f$. Both $\beta_f$ and $\sigma_f$ belong to this class, since their embeddings are into products of locally compact abelian groups and local quasi-convexity is inherited by subgroups. The general compatible-topology and Mackey group framework is treated in \cite{Gabriyelyan2022}.

\begin{proposition}[Least compatible topologies and additive duals]\label{prop:compatible}
The topology $\beta_f$ is least among all Hausdorff group topologies with character group $H_f$. The topology $\sigma_f$ is least among those additionally containing $\tu$. For any compatible Hausdorff TVG topology $\eta$, every continuous real-valued additive map is $x\mapsto cx$. Its real-valued additive dual is therefore either $\{0\}$ or $\{c\,\mathrm{id}:c\in\R\}$. The latter case occurs exactly when $\tu\subseteq\eta$.
\end{proposition}
\begin{proof}
A topology with character group $H_f$ makes all its characters continuous and therefore contains their initial topology $\beta_f$. If it contains $\tu$ as well, it contains $\tu\vee\beta_f=\sigma_f$.

For the additive-dual assertion let $T:(\R,\eta)\to\R_{\mathrm{usual}}$ be continuous and additive. All $p(rT)$ are continuous for $r\in\R$. The proof of Theorem~\ref{thm:rigidity}, using only that its source character set is $H_f$, gives $T=c\,\mathrm{id}$. If $c\ne0$, continuity gives continuity of the identity into the usual real line, hence $\tu\subseteq\eta$. Conversely this inclusion makes every scalar multiple of the identity real-valued continuous.
\end{proof}

The least-member statement does not make $\sigma_f$ finest. To construct a strictly finer compatible TVG we use the following established result: there is a Hausdorff locally quasi-convex topology $\lambda$ on $\Q/\Z$, strictly finer than its circle subspace topology, with exactly the same continuous characters, namely the integer multiples of the inclusion into $\T$ \cite{Barrera2014}. Fix one such $\lambda$. It is not precompact, since a precompact topology is uniquely determined by its continuous character subgroup.

\begin{lemma}[A rational target with unchanged characters]\label{lem:rational-target}
Let $\gamma$ be the initial topology on $\Q$ of
\[
t\longmapsto p(qt)\in(\Q/\Z,\lambda)\qquad(q\in\Q).
\]
Then $\gamma$ is Hausdorff, locally quasi-convex, and invariant under multiplication by every fixed rational scalar. It is not precompact, and
\[
\widehat{(\Q,\gamma)}
=\{t\mapsto p(rt):r\in\Q\}.
\]
Write $\gamma_0$ for the precompact topology on $\Q$ generated by the same maps with their usual circle targets. Then $\gamma_0\subsetneq\gamma$.
\end{lemma}
\begin{proof}
Initial group topologies preserve the stated local quasi-convexity, and rational rescaling permutes the defining maps. Their kernels have zero intersection: for $t\ne0$, choose $q\in\Q$ with $qt=1/2$. Thus the topology is Hausdorff.

A character is continuous if it sends some zero-neighborhood into $[-\frac14,\frac14]+\Z$. For completeness, this usual criterion follows by intersecting inverse images of that neighborhood under multiplication by $1,\ldots,N$; if the first $N$ multiples of a circle value lie in the quarter arc, that value lies within $\frac1{4N}$ of zero. Character continuity for an initial group topology is therefore controlled by finitely many defining maps.

Write their parameters as $q_i=\frac{k_i}d$, with $k_i\in\Z$ and $d\ge1$. All these maps factor continuously through $t\mapsto p(\frac{t}{d})$, since integer multiplication is continuous for $\lambda$. The character is consequently continuous for this single pullback topology. Every element of the kernel $d\Z$ belongs to the closure of zero in this pullback topology: the defining map cannot distinguish it from zero. A continuous character sends that closure into the closure of zero in the Hausdorff group $\T$, which is $\{0\}$. Thus the character vanishes on $d\Z$, factors continuously through $(\Q/\Z,\lambda)$, and has the form $p(\frac{kt}d)$ for some $k\in\Z$. Conversely all rational-parameter characters are among the defining characters. This proves the dual formula.

The map $p:(\Q,\gamma)\to(\Q/\Z,\lambda)$ is continuous and onto. If $\gamma$ were precompact, its image would be precompact, contrary to the choice of $\lambda$. Since $\gamma_0$ is precompact and its defining circle maps are $\gamma$-continuous, strictness follows.
\end{proof}

\begin{lemma}[Finite coefficient quotients]\label{lem:finite-quotient}
Let $h_1,\ldots,h_m\in D$ be rationally independent. Let $\nu$ be the topology on $\R$ generated by the usual coordinate and these maps into $(\Q,\gamma)$. Then
\[
h:(\R,\nu)\longrightarrow(\Q,\gamma)^m,\qquad
h(x)=(h_1(x),\ldots,h_m(x)),
\]
is an open continuous surjection. Every $\nu$-continuous character has the form
\[
x\longmapsto p\left(cx+\sum_{j=1}^m r_jh_j(x)\right),
\qquad c\in\R,\quad r_j\in\Q.
\]
The open-surjection assertion also holds with $\gamma_0$ in place of $\gamma$.
\end{lemma}
\begin{proof}
Independence gives surjectivity by the finite-dimensional annihilator argument. The kernel $N=\bigcap_j\ker h_j$ has finite rational codimension in $\R$, is nonzero, and is usual-dense. For every product neighborhood $W$ of zero in the target and every $\delta>0$,
\[
h\bigl((-\delta,\delta)\cap h^{-1}(W)\bigr)=W:
\]
each fiber is a translate of the dense subgroup $N$ and meets the interval. This proves openness.

The subspace topology of $\nu$ on $N$ is exactly usual. The restriction of a continuous character $\chi$ to $N$ extends to the usual real line, since $N$ is dense and $\T$ is complete. That extension is $p(cx)$ for some $c$. The difference $\chi-p(c\,\cdot)$ vanishes on $N$ and, by the open quotient map $h$, factors as a continuous character on $(\Q,\gamma)^m$. Its restrictions to the factors have rational parameters by Lemma~\ref{lem:rational-target}. Summing those restrictions gives the formula. The same density argument proves openness for $\gamma_0$.
\end{proof}

\begin{theorem}[A strictly finer compatible TVG]\label{thm:finer-compatible}
Let $\eta_f$ be generated by the usual coordinate and all maps
\[
f_a:\R\longrightarrow(\Q,\gamma)\qquad(a\in\R).
\]
Then $\eta_f$ is a Hausdorff nondiscrete locally quasi-convex TVG topology, and
\[
\sigma_f\subsetneq\eta_f,\qquad
\widehat{(\R,\eta_f)}=H_f.
\]
In particular, $\sigma_f$ is not the finest compatible locally quasi-convex TVG topology.
\end{theorem}
\begin{proof}
The initial-topology construction gives local quasi-convexity and the group property. Fixed real dilation permutes the maps $f_a$, proving the TVG property; the usual coordinate makes the topology Hausdorff. Since $\gamma_0\subseteq\gamma$, all defining modular observations are continuous, so $\sigma_f\subseteq\eta_f$ and every character of $H_f$ is continuous.

Conversely, continuity of a character is controlled by the usual coordinate and finitely many coefficient maps, by the quarter-arc criterion in Lemma~\ref{lem:rational-target}. Choose a rational basis $h_1,\ldots,h_m$ of their span. The two finite lists define the same topology with target $\gamma$, because rational linear combinations are continuous for $\gamma$. Each $h_j$ is itself a rational combination of rescalings of $f$, hence is a rescaling of $f$. Lemma~\ref{lem:finite-quotient} gives the character as $p(cx+\sum_jr_jh_j(x))$. The rational sum again combines into one rescaling, so the character belongs to $H_f$. An empty coefficient list gives just an ordinary character. This proves compatibility.

To prove strictness, the map
\[
f:(\R,\sigma_f)\longrightarrow(\Q,\gamma_0)
\]
is an open continuous surjection. Continuity follows from $p(qf(x))=p(f(qx))$. We give the denominator details needed for openness. Let
\[
U=\{x:|x|<\delta,\ \|f_{a_i}(x)\|_{\T}<\varepsilon
\ (1\le i\le m)\}
\]
be a basic $\sigma_f$-neighborhood. Choose a rational basis $h_1,\ldots,h_k$ of the span of $f,f_{a_1},\ldots,f_{a_m}$ with $h_1=f$. Choose a positive integer $d_0$ clearing all denominators of the observation coefficients in this basis, and put $g_j=\frac{h_j}{d_0}$. Then
\[
f_{a_i}=\sum_{j=1}^k m_{ij}g_j,\quad m_{ij}\in\Z,
\qquad f=d_0g_1.
\]
Set
\[
M=1+\sum_{i=1}^m\sum_{j=1}^k|m_{ij}|,
\qquad \varepsilon'=\frac{\varepsilon}{M}.
\]
For any $x$ with $\|g_j(x)\|_{\T}<\varepsilon'$ for every $j$, the integer coefficients give
\[
\|f_{a_i}(x)\|_{\T}
\le\sum_{j=1}^k |m_{ij}|\|g_j(x)\|_{\T}
<\varepsilon.
\]
The maps $g_j$ are rationally independent, so their joint map is onto $\Q^k$. Its kernel has finite rational codimension and is usual-dense. For any $t\in\Q$ with $\|t\|_{\T}<\varepsilon'$, choose a representative with joint coordinates $(t,0,\ldots,0)$ and adjust it by an element of that kernel so that $|x|<\delta$. The resulting point belongs to $U$ and satisfies $f(x)=d_0t$. Consequently
\[
\begin{aligned}
f(U)&\supseteq\{d_0t:t\in\Q,\ \|t\|_{\T}<\varepsilon'\}\\
&=\left\{q\in\Q:\left\|\frac{q}{d_0}\right\|_{\T}<\varepsilon'\right\}.
\end{aligned}
\]
The last set is a $\gamma_0$-neighborhood, because it is the inverse image of a circle arc under the defining character $q\mapsto p(\frac{q}{d_0})$. This proves that images of basic zero-neighborhoods contain target zero-neighborhoods. Translation then shows that the surjection $f$ is open.

If $\eta_f=\sigma_f$, continuity of $f:(\R,\eta_f)\to(\Q,\gamma)$ would descend through this quotient and make the identity from $\gamma_0$ to $\gamma$ continuous. That contradicts $\gamma_0\subsetneq\gamma$. Finally, $\eta_f$ is nondiscrete: its full character set $H_f$ has cardinality $\mathfrak c$, whereas the discrete real line has $2^{\mathfrak c}$ algebraic characters.
\end{proof}

\textbf{Remaining question on the finest compatible topology.}
Does $\mathcal C_f$ have a greatest member? Theorem~\ref{thm:finer-compatible} rules out $\sigma_f$ as that member. The existence or identification of a greatest member is not proved here. This remains a Mackey question restricted to scalar-invariant topologies; the unrestricted locally quasi-convex group problem allows additional topologies. Any proposed supremum still requires a check that it introduces no extra continuous characters.

\section{Maximal nondiscrete topologies}\label{sec:maximal}
A maximal nondiscrete TVG topology has no strictly finer nondiscrete TVG topology. A greatest nondiscrete TVG topology would contain every nondiscrete TVG topology. We establish common properties of maximal members and then construct a pair of nondiscrete TVGs with discrete join. This rules out a greatest member and proves nonuniqueness of maximal extensions.

\begin{lemma}[Countable neighborhoods]
A Hausdorff TVG topology on $\R$ with a countable neighborhood of zero is discrete.
\end{lemma}
\begin{proof}
Choose a countable open set containing zero inside that neighborhood. Its generated additive subgroup $E$ is countable and open. Choose $a\ne0$ outside the countable set of ratios $\frac{h}{k}$ with $h,k\in E$ and $k\ne0$. Then $E\cap aE=\{0\}$. Both subgroups are open because fixed nonzero dilations are homeomorphisms. Thus $\{0\}$ is open.
\end{proof}

\begin{lemma}[Two nondiscrete joins]\label{lem:joins}
Let $\eta$ be a Hausdorff nondiscrete TVG topology on $\R$. Then $\eta\vee\tu$ is nondiscrete. If $\pi$ is any precompact group topology, $\eta\vee\pi$ is also nondiscrete.
\end{lemma}
\begin{proof}
If the first join were discrete, there would be an $\eta$-neighborhood $U$ and $\delta>0$ such that $U\cap(-\delta,\delta)=\{0\}$. Choose an $\eta$-neighborhood $W$ with $W-W\subseteq U$. Distinct members of $W$ have usual distance at least $\delta$, so $W$ is countable. The preceding lemma would make $\eta$ discrete.

For the second assertion, discreteness would give neighborhoods $U$ for $\eta$ and $V$ for $\pi$ with $U\cap V=\{0\}$. Choose neighborhoods $W,Z$ with $W-W\subseteq U$ and $Z-Z\subseteq V$ in the respective topologies. Precompactness covers $\R$ by finitely many translates of $Z$. Each meets $W$ in at most one point, because a difference of two points would lie in $U\cap V$. Thus $W$ is finite, again forcing discreteness.
\end{proof}
The precompact part is a standard obstruction to transversal group topologies; see \cite{DTY2005}. Its elementary proof is included to isolate exactly what is needed for TVGs.

\begin{theorem}[Existence and common properties of maximal members]\label{thm:maximal}
Every Hausdorff nondiscrete TVG topology on $\R$ has a maximal nondiscrete TVG extension. Every maximal member $\mu$ contains both $\tu$ and
\[
B=\theta_{\mathcal A},
\]
the Bohr topology of the discrete additive real line. All algebraic characters are $\mu$-continuous. Every $\mu$-compact subset is finite, every convergent sequence is eventually constant, and $\mu$ is not first countable.
\end{theorem}
\begin{proof}
Order the Hausdorff nondiscrete TVG extensions of a given topology by inclusion. The join of a chain is again a TVG topology. If it were discrete, finitely many zero-neighborhoods from chain members would intersect in $\{0\}$. A single chain member contains these finitely many open sets, contradicting its nondiscreteness. Thus chains have upper bounds in the ordered set, and Zorn's lemma gives a maximal member.

The topology $B$ is a Hausdorff precompact TVG by Theorem~\ref{thm:precompact}. Lemma~\ref{lem:joins} and maximality give $\mu\vee\tu=\mu$ and $\mu\vee B=\mu$. Hence all algebraic characters are continuous.

Apply Glicksberg's theorem to the discrete locally compact group $\R_d$. Its Bohr topology is $B$, so its compact subsets are exactly the finite subsets. The identity from $(\R,\mu)$ into $(\R,B)$ is continuous. Every $\mu$-compact subset is therefore finite. A convergent sequence and its limit form a compact set, so Hausdorffness forces the sequence to be eventually equal to its limit.

Finally, a first countable nondiscrete Hausdorff group admits a sequence of nonzero points tending to zero: use a nested countable neighborhood base and choose a nonzero point in each member. This contradicts the preceding conclusion. Thus $\mu$ is not first countable.
\end{proof}
The topology $B$ differs from the usual-line Bohr topology $t$: $B$ uses every algebraic character, whereas $t$ uses only usual-continuous characters. In particular, the sequence $1/n$ converges in $t$ but not in $B$.

\begin{lemma}[Disjoint rational subspaces with all finite scalar patterns]\label{lem:patterns}
There are rational subspaces $E,F\subset\R$ with $E\cap F=\{0\}$ such that for every finite nonempty $A\subset\R\setminus\{0\}$ there are nonzero $u,v\in\R$ satisfying
\[
Au\subset E,\qquad Av\subset F.
\]
\end{lemma}
\begin{proof}
Enumerate the finite nonempty subsets of $\R\setminus\{0\}$ in $\mathfrak c$ stages, using the initial ordinal of that cardinal. Build increasing rational subspaces $E,F$ by transfinite recursion, keeping them disjoint. At the stage for $A$, let $L=E+F$ be the span of all previously selected vectors and put $V=\operatorname{span}_{\Q}A$. The space $V$ is finite-dimensional over $\Q$ and hence countable. Before this stage, fewer than $\mathfrak c$ vectors have been selected, so $|L|<\mathfrak c$.

The forbidden set
\[
\bigcup_{w\in V\setminus\{0\}}w^{-1}L
\]
has cardinality at most $\max\{\aleph_0,|L|\}<\mathfrak c$. Choose $u$ outside it. Then $uV\cap L=\{0\}$, so adding $uV$ to $E$ preserves disjointness from $F$. Apply the same choice with the enlarged span to obtain $v$, and add $vV$ to $F$. At limit stages take unions. Each addition has finite rational dimension, so the cardinal estimate remains valid at every stage; regularity of $\mathfrak c$ is not assumed. The resulting subspaces are disjoint and contain the selected copies of every $A$.
\end{proof}

\begin{theorem}[No greatest nondiscrete TVG]\label{thm:no-greatest}
There are Hausdorff nondiscrete TVG topologies $\rho_E,\rho_F$, each strictly finer than $\tu$, whose join is discrete. Consequently no greatest nondiscrete Hausdorff TVG topology exists on $\R$, even within the class of topologies finer than usual.
\end{theorem}
\begin{proof}
Choose $E,F$ as in Lemma~\ref{lem:patterns}. Define a zero-neighborhood base for $\rho_E$ by
\[
U_E(\delta,A)=(-\delta,\delta)\cap
\bigcap_{a\in A}a^{-1}E,
\]
where $\delta>0$ and $A\subset\R\setminus\{0\}$ is finite, possibly empty. Define $\rho_F$ similarly. The subspace conditions are subgroup conditions, so halving the interval bound gives a neighborhood whose sum with itself lies in the original neighborhood. Inverses preserve each condition. Translates give the corresponding group topologies.

Multiplication by a fixed $c\ne0$ changes the interval bound to $\frac{\delta}{|c|}$ and changes $A$ to $cA$. Thus every fixed dilation is continuous, and these are TVG topologies. They contain $\tu$ and are Hausdorff.

For nonempty $A$, choose $u\ne0$ with $Au\subset E$. Then the rational line $u\Q$ lies in $\bigcap_{a\in A}a^{-1}E$. It has arbitrarily small nonzero elements, so every $U_E(\delta,A)$ contains a nonzero point. The empty case is immediate. Hence $\rho_E$ is nondiscrete, and the same proof applies to $\rho_F$.

The subspaces $E,F$ are open in their respective topologies: intersect each with a sufficiently large interval to obtain a neighborhood at each of its points. Both are nonzero proper rational subspaces, since they are disjoint and each contains a selected nonzero pattern. A proper additive subgroup cannot be usual-open. Thus both topologies are strictly finer than usual.

Their join makes $E\cap F=\{0\}$ open and is discrete. A greatest nondiscrete topology would contain this join, which is impossible.
\end{proof}

\begin{corollary}[Nonuniqueness of maximal nondiscrete extensions]\label{cor:nonunique-maximal}
There are distinct maximal nondiscrete Hausdorff TVG topologies on $\R$. In particular, the common nondiscrete topology $\tu\vee B$ has at least two different maximal nondiscrete extensions.
\end{corollary}
\begin{proof}
Extend $\rho_E$ and $\rho_F$ separately by Theorem~\ref{thm:maximal}. A common nondiscrete extension would contain their discrete join, so these maximal extensions are distinct. Both contain $\tu$ and $B$ by that same theorem. The topology $\tu\vee B$ is nondiscrete by Lemma~\ref{lem:joins}, since $B$ is precompact.
\end{proof}
This proves nonuniqueness in general, including for the specified common starting topology. It does not assert that every prescribed nondiscrete TVG has multiple maximal extensions. All maximal members still have the same full algebraic character group $\mathcal A$, as shown above.

\section{Compatible joins and maximal compatible extensions}\label{sec:joins-v8}
We retain $H_f$ and $\mathcal C_f$ from Section~\ref{sec:compatible}. A greatest member contains every member of the class; a maximal member merely has no strictly finer member in the class. The next argument adapts the standard compatible-topology framework \cite{Gabriyelyan2022,AD2021} to fixed-scalar TVGs.

\begin{proposition}[Directed joins and the greatest-member criterion]\label{prop:directed-v8}
Every nonempty directed join in $\mathcal C_f$ belongs to $\mathcal C_f$. Every member has a maximal extension in $\mathcal C_f$. The following conditions are equivalent: existence of a greatest member; compatibility of every binary join; uniqueness of the maximal member.
\end{proposition}
\begin{proof}
Joins preserve group topologies, local quasi-convexity, Hausdorffness, and continuity of fixed dilations, by the diagonal embedding into a product. All characters of $H_f$ stay continuous. A character continuous for a directed join maps some basic zero-neighborhood into the closed quarter arc. This neighborhood uses finitely many family members, dominated by one member. The quarter-arc criterion of Lemma~\ref{lem:rational-target} makes the character continuous for that member, hence an element of $H_f$.

Zorn's lemma now gives maximal extensions, since chains are directed. A greatest member bounds every binary join, so no join can acquire an extra character. Conversely, compatibility of binary joins gives compatibility of all finite joins; their directed supremum is compatible and greatest. Finally a unique maximal member contains every member, since each extends to a maximal member. A greatest member is plainly the unique maximal member.
\end{proof}
Checking a single selected pair cannot prove compatibility of all binary joins.

\section{Rational quotients and transfer to the real line}\label{sec:rational-v8}
Define
\[
\chi_r(q)=p(rq),\quad
\Gamma=\{\chi_r:r\in\Q\},\qquad
\gamma_0=\operatorname{init}\{\chi_r:r\in\Q\}.
\]
The topology $\gamma_0$ is Hausdorff, precompact and metrizable, with full character set $\Gamma$, by precompact duality. Let $\mathcal L_\Gamma$ be the Hausdorff locally quasi-convex group topologies on $\Q$ with full character set $\Gamma$. Let $\mathcal Q_\Gamma$ be those members admitting every fixed rational multiplication continuously. The topology $\gamma$ of Lemma~\ref{lem:rational-target} belongs to $\mathcal Q_\Gamma$.

\begin{lemma}[The quotient associated with a coefficient]\label{lem:quotient-v8}
For $\eta\in\mathcal C_f$, $E=\ker f$ is closed. Identify $\R/E$ with $\Q$ through $f$. The Hausdorff quotient topology has exactly the characters $\Gamma$ and admits all fixed rational multiplications continuously. Local quasi-convexity of this quotient is not asserted.
\end{lemma}
\begin{proof}
We have $E=\bigcap_{r\in\Q}\ker(p\circ rf)$: if $f(x)\ne0$, choose rational $r$ with $rf(x)=\frac{1}{2}$. These are characters in $H_f$, so $E$ is closed.

A quotient character corresponds to $\chi_{c,a}$ vanishing on $E$. Then $h(x)=cx+f(ax)$ sends the divisible group $E$ into $\Z$, hence is zero on $E$. Thus $cE\subseteq\Q$. Since $E$ is uncountable, $c=0$. The remaining condition $f_a(E)=0$ implies $f_a=rf$, where $r=f_a(u)\in\Q$ for $f(u)=1$. Conversely each $p(rf)$ vanishes on $E$. This identifies the quotient dual. Rational dilations descend because $E$ is a rational subspace.
\end{proof}

\begin{theorem}[Rational transfer and incompatible joins]\label{thm:transfer-v8}
For $\delta\in\mathcal Q_\Gamma$, let $T_f(\delta)$ be generated by the usual coordinate and all maps $f_a:\R\to(\Q,\delta)$. Then $T_f(\delta)\in\mathcal C_f$, and $f:T_f(\delta)\to(\Q,\delta)$ is an open surjection. This assignment is an order embedding, and
\[
T_f(\delta_1)\vee T_f(\delta_2)
=T_f(\delta_1\vee\delta_2).
\]
The transferred join is compatible if and only if the rational join has character set $\Gamma$. An incompatible rational pair would rule out a greatest member of $\mathcal C_f$.
\end{theorem}
\begin{proof}
Lemma~\ref{lem:finite-quotient} applies with $\delta$ replacing $\gamma$: its proof uses only rational scalar continuity, the dual $\Gamma$, and density of kernels of finitely many independent rational functionals. Finite rational coordinate changes are continuous for $\delta$. The character argument of Theorem~\ref{thm:finer-compatible} therefore proves compatibility; the initial construction gives the other structural properties.

For openness, choose a rational basis $h_1=f,h_2,\ldots,h_m$ containing any finite list of observations in its span. The joint map is onto $\Q^m$ with usual-dense kernel. The image under $f$ of a basic source neighborhood contains the first-coordinate projection of a target neighborhood in $(\Q,\delta)^m$: adjust each fiber into the prescribed usual interval. That projection is a $\delta$-neighborhood. Thus $f$ is open and recovers $\delta$ as quotient topology, proving order reflection as well as preservation.

The equality of joins follows from their defining maps. If the rational join is compatible, the same finite-coordinate proof applies. Otherwise take an extra character $\psi$ on that join. Then $\psi\circ f$ is continuous on the transferred join. The annihilator calculation of Lemma~\ref{lem:quotient-v8} shows that $\psi\circ f\in H_f$ would imply $\psi\in\Gamma$, a contradiction. A greatest compatible member cannot dominate an incompatible join.
\end{proof}

\begin{proposition}[Cofinal rational scalar continuity]\label{prop:rationalize-v8}
Every $\delta\in\mathcal L_\Gamma$ has a finer member $\delta^\Q\in\mathcal Q_\Gamma$. Consequently a greatest member of $\mathcal Q_\Gamma$ exists exactly when $(\Q,\Gamma)$ admits a Mackey topology in the full locally quasi-convex group class.
\end{proposition}
\begin{proof}
Let $\delta_n$ be the pullback of $\delta$ by $q\mapsto \frac{q}{n}$, and put $\delta^\Q=\bigvee_{n\ge1}\delta_n$. A character $\psi$ is $\delta_n$-continuous exactly when $n\psi\in\Gamma$. The group $\operatorname{Hom}(\Q,\T)$ is torsion-free: $m\psi=0$ implies $\psi(q)=m\psi(\frac{q}{m})=0$. Thus $n\psi=\chi_r$ forces $\psi=\chi_{\frac{r}{n}}$, proving compatibility of $\delta_n$.

If $d$ is a multiple of $n$, integer multiplication gives $\delta_n\subseteq\delta_d$. The family is directed, and its join remains compatible by the quarter-arc argument of Proposition~\ref{prop:directed-v8}. Multiplication by $\frac{a}{b}$, followed by $q\mapsto \frac{q}{n}$, is the integer multiple $a$ of the defining map $q\mapsto \frac{q}{bn}$. Hence all rational multiplications are continuous. Also $\delta=\delta_1\subseteq\delta^\Q$. Cofinality proves the final equivalence; a greatest member of the larger class equals its own rationalized refinement.
\end{proof}

\begin{corollary}[A necessary rational Mackey condition]\label{cor:necessary-v8}
If $\mathcal C_f$ has a greatest member, $(\Q,\Gamma)$ admits a Mackey topology.
\end{corollary}
\begin{proof}
Let $\mu$ be greatest and $\nu$ its quotient topology through $f$. Theorem~\ref{thm:transfer-v8} makes $\nu$ dominate every $\delta\in\mathcal Q_\Gamma$; Lemma~\ref{lem:quotient-v8} gives its dual $\Gamma$. Take the locally quasi-convex modification $\nu_{\mathrm{lqc}}$, the finest locally quasi-convex topology coarser than $\nu$ \cite{AD2021}. Character preimages of small circle arcs are quasi-convex neighborhoods, so this modification preserves the characters. It is Hausdorff since $\Gamma$ separates points. Rational scalar continuity persists because inverse images of quasi-convex neighborhoods under continuous homomorphisms are quasi-convex neighborhoods. The modification dominates every locally quasi-convex $\delta\subseteq\nu$. It is therefore greatest in $\mathcal Q_\Gamma$, and Proposition~\ref{prop:rationalize-v8} applies.
\end{proof}
The converse is not proved. The usual-topology conclusion in \cite{Barrera2014} does not assert that no Mackey topology exists for $(\Q,\Gamma)$.

\section{Ultrafilter tests for compatible joins}\label{sec:ultrafilters-v8}
A family of characters is equicontinuous if, for every circle zero-neighborhood $V$, one source zero-neighborhood $U$ satisfies $\chi(U)\subseteq V$ for every character in the family.

An ultrafilter $\mathcal U$ on a set $I$ is a family of subsets closed under finite intersections and supersets, excluding the empty set, such that exactly one of $A$ and $I\setminus A$ belongs to it for every $A\subseteq I$. It is free if it contains no finite set. The notation $\lim_{i\to\mathcal U}z_i=z$ means that $\{i:z_i\in V\}\in\mathcal U$ for every neighborhood $V$ of $z$. Every family in a compact Hausdorff space has a unique ultrafilter limit. Extending a family with the finite intersection property uses the ultrafilter lemma.

\begin{proposition}[Ultralimits and a sufficient counterexample test]\label{prop:ultralimit-v8}
For $r_i\in\Q$, the pointwise limit $\psi(q)=\lim_{i\to\mathcal U}p(r_iq)$ is an algebraic character. Suppose $\delta_1,\delta_2\in\mathcal Q_\Gamma$, $\psi\notin\Gamma$, and respective zero-neighborhoods satisfy
\[
\|r_iq\|_{\T}\le\frac14
\quad(q\in U_1\cap U_2,\ i\in I).
\]
Then the join is incompatible. If $r_n\in\Q$ tends usually to an irrational $c$, every free ultrafilter gives $\psi(q)=p(cq)\notin\Gamma$.
\end{proposition}
\begin{proof}
Compactness supplies the limits and circle addition preserves additivity. The closed quarter-arc bound passes to $\psi$, whose continuity in the join follows from the quarter-arc criterion. Usual convergence gives the stated limit. Equality $p(cq)=p(rq)$ on all of $\Q$ would make $(c-r)\Q$ a divisible subgroup of $\Z$, hence zero, contrary to irrationality.
\end{proof}
This supplies a sufficient test, not the required neighborhoods. No infinite subsequence of the characters associated with $r_n\to c\notin\Q$ can be equicontinuous for either compatible topology: its limit would already be continuous there. Splitting those approximants into individually equicontinuous infinite families cannot solve the construction.

\begin{theorem}[An exact ultrafilter criterion]\label{thm:ultrafilter-v8}
For an algebraic character $\psi$ put
\[
A_\psi=\{q\in\Q:\|\psi(q)\|_{\T}>\tfrac14\},
\]
where the norm on a circle element is distance from zero. A join of $\delta_1,\delta_2\in\mathcal Q_\Gamma$ is compatible exactly when every $\psi\notin\Gamma$ admits an ultrafilter $\mathcal F$ on $\Q$ converging to zero in both topologies and containing $A_\psi$.
\end{theorem}
\begin{proof}
If the join is compatible, every joint zero-neighborhood meets $A_\psi$, by discontinuity and the quarter-arc criterion. The sets $U_1\cap U_2\cap A_\psi$ have the finite intersection property. An extending ultrafilter contains every zero-neighborhood from both topologies and contains $A_\psi$. Conversely such an ultrafilter converges to zero in the join. Its membership of $A_\psi$ contradicts continuity of $\psi$, by the inverse image of the open quarter arc.
\end{proof}
The proof applies equally to $\mathcal C_f$, replacing $\Q,\Gamma$ by $\R,H_f$. No selective ultrafilter or continuum hypothesis is used.

\section{Restrictions on a rational counterexample}\label{sec:restrictions-v8}
A group topology is called locally precompact here when its group completion is locally compact.

\begin{theorem}[Reduction to metrizable topologies]\label{thm:metric-pair-v8}
If two members of $\mathcal Q_\Gamma$ have an incompatible join, two metrizable members do, witnessed by the same extra character.
\end{theorem}
\begin{proof}
Choose $\psi\notin\Gamma$ continuous for $\delta_1\vee\delta_2$ and respective neighborhoods $U_i$ whose intersection maps into the quarter arc. Enumerate rational scalars $s_1,s_2,\ldots$, and choose a countable zero-neighborhood base $(B_n)$ for $\gamma_0$. Recursively choose symmetric quasi-convex $\delta_i$-neighborhoods $W_{i,n}$, starting with $W_{i,0}\subseteq U_i$, such that
\[
\begin{aligned}
W_{i,n}+W_{i,n}&\subseteq W_{i,n-1},\\
W_{i,n}&\subseteq W_{i,n-1}\cap B_n,\\
s_jW_{i,n}&\subseteq W_{i,n-1}\quad(1\le j\le n).
\end{aligned}
\]
Continuity and local quasi-convexity permit each choice. These sets define a first countable Hausdorff group topology $\zeta_i$ between $\gamma_0$ and $\delta_i$. Its dual is $\Gamma$, so its defining sets remain quasi-convex. The last inclusion gives rational scalar continuity. The Birkhoff--Kakutani theorem gives metrizability. The retained $W_{i,0}$ keep $\psi$ continuous in the new join by the quarter-arc criterion.
\end{proof}

\begin{proposition}[No locally precompact participant]\label{prop:local-precompact-v8}
The only member of $\mathcal L_\Gamma$ with locally compact completion is $\gamma_0$. Neither participant of an incompatible rational pair can be locally precompact.
\end{proposition}
\begin{proof}
For such a completion $L$, extension of characters gives countable $\widehat L$. This locally compact Hausdorff group is Baire and must be discrete: otherwise it is a countable union of nowhere dense singletons. Pontryagin duality makes $L$ compact. The original topology is precompact and is uniquely determined by its character group, so it is $\gamma_0$. Joining this least compatible topology to any other compatible topology changes nothing.
\end{proof}

\begin{lemma}[Torsion and roots cannot create an extra character]\label{lem:roots-v8}
The group $\operatorname{Hom}(\Q,\T)$ is torsion-free. If $m\psi\in\Gamma$ for a positive integer $m$, then $\psi\in\Gamma$. For characters on $\R$, $m\psi\in H_f$ similarly implies $\psi\in H_f$.
\end{lemma}
\begin{proof}
For a divisible group $G$, $m\psi=0$ implies $\psi(x)=m\psi(\frac{x}{m})=0$. If $m\psi=\chi_r$, subtract $\chi_{\frac{r}{m}}$. For $H_f$, subtract $\chi_{\frac{c}{m},\frac{a}{m}}$ when $m\psi=\chi_{c,a}$ and use the same torsion-free argument.
\end{proof}
Thus an extra character must give a new direction modulo the prescribed dual. If the rational topologies are coarser than usual, an extra character has the form $p(cq)$ with irrational $c$. Its rational multiples and $\Gamma$ form the direct sum $\Gamma+\Q\chi_c$; this is still countable, giving no cardinality contradiction.

Gabriyelyan's counterexample \cite{Gabriyelyan2018} uses a different mechanism. On the Graev free abelian group of a convergent sequence, algebraically $\Z^{(\mathbb N)}$, two compatible refinements have a join making a nonzero finite-order character continuous. The order-two twist gives the character
\[
(k_n)\longmapsto p\left(\frac12\sum_n k_n\right).
\]
Their combined neighborhood conditions force every $k_n$ to be even. It is discontinuous originally: the coordinate vectors tend to zero while its value on each is $\frac{1}{2}+\Z$. This is the mechanism in Claims 1--2 of \cite{Gabriyelyan2018}. Lemma~\ref{lem:roots-v8} prevents direct transfer of that finite-order step to $\Q$ or $\R$.

\section{An explicit pair from prime square denominators}\label{sec:prime-v8}
For $p=2$ or $3$, $d\in\mathbb N_{>0}$ and $q\in\Q$, put
\[
P_{p,d}(q)=\sup_{n\ge0}
\left\|\frac{q}{d\,p^{n^2}}\right\|_{\T}.
\]
Here the letter $p$ in the subscript denotes the prime, whereas $p(t)=t+\Z$ elsewhere denotes the quotient map; their arguments distinguish the uses. To avoid ambiguity in the following proof, characters are also written directly as $cq+\Z$. Let $\gamma_p$ be generated by the functions $P_{p,d}$, with basic zero-neighborhoods
\[
U_p(d,\varepsilon)=\{q:P_{p,d}(q)<\varepsilon\}.
\]
These are invariant metric balls for the corresponding individual metrics. Common multiples of denominators and smaller radii refine every finite intersection.

\begin{theorem}[Compatibility of the two prime constructions]\label{thm:prime-v8}
The topologies $\gamma_2,\gamma_3$ belong to $\mathcal Q_\Gamma$. They are metrizable, not precompact, and not locally precompact.
\end{theorem}
\begin{proof}
Triangle inequalities on the circle give the group topology. Each individual $P_{p,d}$ separates points, since for nonzero $q$ its terms $\frac{q}{d p^{n^2}}$ are eventually nonzero and smaller than $\frac{1}{2}$ in absolute value. Countably many denominators and radii give metrizability. Uniform circle observations give local quasi-convexity; equivalently use the closed balls of radii $\frac{1}{4m}$, intersections of character preimages of quasi-convex arcs. For integers $a$ and positive integers $b,d$,
\[
P_{p,d}\left(\frac{a}{b}q\right)
\le |a|P_{p,bd}(q),
\]
so fixed rational multiplications are continuous.

The character computation uses the following external theorem. For a sequence of positive integers $b_0=1$ with $b_n\mid b_{n+1}$, strictly increasing and with $\frac{b_{n+1}}{b_n}\to\infty$, the topology on $\Z$ of uniform convergence on $k\mapsto \frac{k}{b_n}+\Z$ has character set generated by $\frac{1}{b_n}+\Z$ \cite[Theorem 4.4]{AB2012}. Apply it to $b_n=p^{n^2}$. Its dual is $\Z[1/p]/\Z$, where $\Z[1/p]=\{\frac{m}{p^j}:m\in\Z,\ j\ge0\}$.

A $\gamma_p$-continuous character is controlled by one $P_{p,d}$, by a common denominator and the quarter-arc criterion. Since $P_{p,d}(q)\le \frac{|q|}{d}$, it is usual-continuous on $\Q$. Extension to $\R$ therefore makes it $cq+\Z$ for a real $c$. Restricting to $q=dk$ gives the integer uniform topology just described. Thus $cd+\Z\in\Z[1/p]/\Z$, so $c\in\Q$. Conversely $\frac{mq}{d}+\Z$ is continuous, since the $n=0$ term controls $\frac{q}{d}+\Z$. The full dual is exactly $\Gamma$.

The infinite character family $q\mapsto \frac{q}{p^{n^2}}+\Z$ is equicontinuous, controlled by $P_{p,1}$. A precompact group has only finite equicontinuous character sets: extend to its compact completion, whose discrete dual has compact, hence finite, neighborhood polars. Thus $\gamma_p$ is not precompact. Proposition~\ref{prop:local-precompact-v8} excludes local precompactness.
\end{proof}

\begin{theorem}[Incomparability]\label{thm:prime-incomparable-v8}
Neither $\gamma_2\subseteq\gamma_3$ nor $\gamma_3\subseteq\gamma_2$ holds.
\end{theorem}
\begin{proof}
Fix $U_3(d,\varepsilon)$ and write $d=2^v e$ with $e$ odd. Choose $n$ with $n^2>v$, put $m=2^{n^2-1-v}$, and set $q_k=dm3^{k^2}$. Terms with index at most $k$ are integers in $P_{3,d}(q_k)$; later terms have absolute value at most $m3^{-(2k+1)}$. Hence
\[
P_{3,d}(q_k)\le m3^{-(2k+1)}\longrightarrow0.
\]
But $\frac{q_k}{2^{n^2}}=\frac{e3^{k^2}}{2}$ is a half-integer, so $P_{2,1}(q_k)\ge\frac{1}{2}$. Thus no $\gamma_3$-neighborhood is contained in $U_2(1,\frac{1}{4})$, proving $\gamma_2\not\subseteq\gamma_3$.

For the reverse direction write $d=3^v e$ with $3$ not dividing $e$, choose $n^2>v$, put $m=3^{n^2-1-v}$, and use $q_k=dm2^{k^2}$. Then $P_{2,d}(q_k)\le m2^{-(2k+1)}\to0$, whereas $\frac{q_k}{3^{n^2}}=\frac{e2^{k^2}}{3}$ has circle norm $\frac{1}{3}$. This proves the other noninclusion.
\end{proof}

\begin{proposition}[The join is still weaker than usual]\label{prop:prime-join-v8}
Every zero-neighborhood of $\gamma_2\vee\gamma_3$ contains arbitrarily large positive integers. Consequently the join is strictly coarser than the usual topology on $\Q$, and its restriction to $\Z$ is nondiscrete.
\end{proposition}
\begin{proof}
The bound $P_{p,d}(q)\le \frac{|q|}{d}$ gives coarseness. A common denominator reduces a joint neighborhood to
\[
P_{2,d}(q)<\varepsilon,\qquad P_{3,d}(q)<\varepsilon.
\]
Set $q=dz$. Given a positive integer $M$, choose an integer $L>2$ with $\frac{1}{L}<\varepsilon$. Simultaneous Dirichlet approximation, by partitioning a finite-dimensional unit cube into $L$ pieces in each coordinate, gives an integer
\[
1\le z\le L^{2(k+1)+M}
\]
whose multiples of all the following numbers have circle norm at most $\frac{1}{L}$:
\[
2^{-n^2},\ 3^{-n^2}\ (0\le n\le k),
\qquad 2^{-j}\ (1\le j\le M).
\]
If $z\le M$, let $v_2(z)$ denote the largest nonnegative integer $r$ for which $2^r$ divides $z$, and put $j=v_2(z)+1$. Since $1\le j\le 2^{v_2(z)}\le z\le M$, the auxiliary coordinate $2^{-j}$ is present. The integer $z/2^{v_2(z)}$ is odd, so $\frac{z}{2^j}\equiv\frac12\pmod1$ and $\left\|\frac{z}{2^j}\right\|_{\T}=\frac12$, contradicting the bound $\frac1L<\frac12$. Thus $z>M$. For $n\ge k+1$ the remaining observations have norm at most
\[
\frac{L^{2(k+1)+M}}{2^{(k+1)^2}},
\]
which is smaller than $\varepsilon$ for sufficiently large $k$. Thus $q=dz$ lies in the joint neighborhood and exceeds $dM$. Unboundedness proves strictness, because the usual topology has bounded zero-neighborhoods.
\end{proof}
Incomparability and strict coarseness do not compute the dual of the join. A strictly weaker topology than usual may still have irrational-parameter characters.

\begin{proposition}[Reduction of the candidate to an integer polar problem]\label{prop:integer-v8}
On $\Z$ put $P(k)=\max\{P_{2,1}(k),P_{3,1}(k)\}$ and let $\omega$ be its metric topology. Define
\[
H=\{\alpha+\Z\in\T:k\mapsto\alpha k+\Z
\text{ is }\omega\text{-continuous}\}.
\]
Then $\gamma_2\vee\gamma_3$ is compatible with $\Gamma$ if and only if $H\subseteq\Q/\Z$.
\end{proposition}
\begin{proof}
Every character of the rational join is usual-continuous, hence $cq+\Z$. A common denominator shows that its continuity is equivalent to continuity of $t\mapsto cd\,t+\Z$ on $\Q$ with the metric $P(t)=\max\{P_{2,1}(t),P_{3,1}(t)\}$, for some positive integer $d$.

For a real $\alpha$, continuity of $\alpha t+\Z$ for this rational metric is equivalent to continuity of its restriction to $\Z$. One direction is immediate. For the other, if $P(t)<\delta<\frac{1}{2}$, the $n=0$ observation gives an integer $k$ with $|t-k|<\delta$. Since $P(t-k)\le |t-k|$, we have $P(k)<2\delta$. Consequently
\[
\|\alpha t\|_{\T}
\le \|\alpha k\|_{\T}+|\alpha|\delta.
\]
Continuity on the integer subgroup therefore proves continuity on $\Q$.

If $H$ is torsion, every continuous rational-join character satisfies $cd+\Z\in\Q/\Z$, forcing $c\in\Q$. Conversely an irrational element $\alpha+\Z\in H$ gives the extra character $\alpha q+\Z$ already for the metric $P$ and hence for the rational join.
\end{proof}
This supplies an exact remaining test: with
\[
B_\varepsilon=\{k\in\Z:P(k)<\varepsilon\},\qquad
B_\varepsilon^\triangleright=
\{\alpha+\Z:\|\alpha k\|_{\T}\le\frac{1}{4}
\text{ for all }k\in B_\varepsilon\},
\]
one must decide whether every such polar consists of torsion elements. The quarter-arc criterion identifies $H$ with the union of these polars. No torsion conclusion for these mixed-prime polars is proved here. The single-sequence theorem \cite[Theorem 4.4]{AB2012} does not by itself cover their union.

\section{Counterexample criteria and unresolved questions}\label{sec:status-v8}
The prime-square topologies have the proved properties in Section~\ref{sec:prime-v8}, but their join dual remains uncomputed. By contrast, an alternating-block attempt already introduces irrational characters into each individual topology and therefore fails compatibility. Appendix~\ref{app:blocks-v11} records that construction and the full proof of its failure; no main construction depends on it.

The next proposition isolates a sufficient witness for the prime pair. Its hypothesis is an unresolved condition, not an assertion that such a witness exists.

\begin{proposition}[What an irrational witness would prove]\label{prop:conditional-v8}
If an irrational $c$ satisfies
\[
\max\{P_{2,d}(q),P_{3,d}(q)\}<\varepsilon
\quad\Longrightarrow\quad
\|cq\|_{\T}\le\frac14
\]
for all $q\in\Q$, with some $d\ge1$ and $\varepsilon>0$, then $\mathcal C_f$ has no greatest member.
\end{proposition}
\begin{proof}
The quarter-arc criterion makes $q\mapsto cq+\Z$ continuous in $\gamma_2\vee\gamma_3$, but it is not in $\Gamma$. Apply Theorem~\ref{thm:transfer-v8}. In particular the transferred join makes $x\mapsto cf(x)+\Z$ continuous. This character cannot lie in $H_f$: the annihilator calculation in Lemma~\ref{lem:quotient-v8} would force its rational quotient character to belong to $\Gamma$.
\end{proof}
The hypothesis of Proposition~\ref{prop:conditional-v8} has not been established. Its implication is a proved conditional statement, not a nonexistence theorem for $\mathcal C_f$. Conversely, proving compatibility of the selected prime pair would not prove that all joins in $\mathcal C_f$ remain compatible. The full greatest-compatible-topology question therefore remains unresolved by this investigation. Maximal compatible extensions exist by Proposition~\ref{prop:directed-v8}, irrespective of that unanswered question. The already proved nonexistence of a greatest nondiscrete TVG in Theorem~\ref{thm:no-greatest} concerns a different class.

The main existence problem can be stated in theorem form while its status remains explicitly unproved.

\begin{conjecture}[Greatest Compatible TVG Conjecture]\label{conj:greatest-v12}
For every nonzero $f\in D$, there exists a topology $\mu_f\in\mathcal C_f$ such that
\[
\eta\subseteq\mu_f\qquad\text{for every }\eta\in\mathcal C_f.
\]
Here $\mathcal C_f$ is the class of Hausdorff locally quasi-convex TVG topologies on $\R$ whose full continuous character group is $H_f$.
\end{conjecture}
This is the \emph{Greatest Compatible TVG Conjecture}, proposed as an open existence assertion. It is not an established theorem, and neither a proof nor a counterexample is supplied here. For each fixed $f$, Proposition~\ref{prop:directed-v8} makes the assertion equivalent to compatibility of every binary join in $\mathcal C_f$. The already proved existence of maximal members does not establish this greatest-member conclusion. The corresponding restricted compact-set question is recorded separately in Appendix~\ref{app:segments-v12}. Section~\ref{sec:dilations-v13} relates the conjecture to the unrestricted Mackey problem for the same dual pair.

\section{Scalar invariance and lifting from the torsion quotient}\label{sec:scalar-v10}
The target class remains $\mathcal C_f$, whose members are Hausdorff locally quasi-convex TVG topologies on $\R$ with full character set $H_f$. Scalar continuity and compatibility are separate requirements.

\begin{proposition}[The property preserved by an arbitrary join]\label{prop:scalar-join-v10}
If $\eta_1,\eta_2\in\mathcal C_f$, then $\eta_1\vee\eta_2$ is a Hausdorff locally quasi-convex TVG topology and every member of $H_f$ remains continuous. The join belongs to $\mathcal C_f$ exactly when it has no additional continuous characters.
\end{proposition}
\begin{proof}
The diagonal map into $(\R,\eta_1)\times(\R,\eta_2)$ induces the join topology. Subgroups and products preserve the Hausdorff group and local quasi-convexity properties. For each fixed $a\in\R$, dilation by $a$ is continuous in both factors and hence in the induced topology. Refining a topology preserves continuity of its characters. Membership in $\mathcal C_f$ therefore reduces precisely to equality of the resulting dual with $H_f$.
\end{proof}
Directed joins have no additional characters by Proposition~\ref{prop:directed-v8}. For a general binary join that conclusion requires proof; TVG invariance alone does not provide it.

Let $A=\Q/\Z$, let $\pi:\Q\to A$ be $\pi(t)=t+\Z$, and let $\iota:A\to\T$ be the inclusion. The circle topology on $A$ is its subspace topology in $\T$. Its continuous characters are exactly $n\iota$ for $n\in\Z$: they extend to its completion $\T$, whose characters are integer multiplications. Write $\mathcal L_A$ for the Hausdorff locally quasi-convex group topologies on $A$ with exactly this character set. Every such topology contains the circle topology, since $\iota$ is continuous.

\begin{lemma}[Lifting a compatible quotient topology]\label{lem:quotient-lift-v10}
For $\lambda\in\mathcal L_A$, let $\gamma_\lambda$ be the initial topology on $\Q$ of the maps
\[
t\longmapsto\pi(rt)\in(A,\lambda)\qquad(r\in\Q).
\]
Then $\gamma_\lambda\in\mathcal Q_\Gamma$.
\end{lemma}
\begin{proof}
Initial topologies of this form are locally quasi-convex group topologies. The kernels separate points: for $t\ne0$, take $r=\frac{1}{2t}$. Composition with a fixed rational multiplication replaces $r$ by another rational parameter, proving scalar continuity.

Let $\psi:\Q\to\T$ be a continuous character. Its quarter-arc bound is controlled by finitely many maps with parameters $r_i=\frac{k_i}{d}$, where $d$ is a positive integer and $k_i\in\Z$. Integer multiplication is continuous on $(A,\lambda)$, so all these maps factor continuously through
\[
h_d:\Q\to(A,\lambda),\qquad h_d(t)=\pi\left(\frac{t}{d}\right).
\]
The same bound makes $\psi$ continuous for the single pullback topology of $h_d$, by the quarter-arc criterion. Its kernel $d\Z$ lies in the closure of zero in that topology. Since $\T$ is Hausdorff, $\psi$ kills $d\Z$. The surjective map $h_d$ is a quotient map for its pullback topology, so $\psi$ descends to a $\lambda$-continuous character of $A$. This character is $n\iota$ for some $n\in\Z$, whence $\psi(t)=p(\frac{nt}{d})\in\Gamma$. Conversely, every $\chi_r\in\Gamma$ is the composition of a defining map with $\iota$ and is continuous. This proves the full dual formula.
\end{proof}
This is the general form of Lemma~\ref{lem:rational-target}; strict refinement or nonprecompactness of $\lambda$ is unnecessary for this lifting assertion.

\begin{theorem}[A sufficient quotient counterexample]\label{thm:quotient-pair-v10}
Suppose $\lambda_1,\lambda_2\in\mathcal L_A$ and the join $\lambda_1\vee\lambda_2$ has a continuous character $\psi$ that is not $n\iota$ for any $n\in\Z$. Then $\gamma_{\lambda_1}\vee\gamma_{\lambda_2}$ has a character outside $\Gamma$. Consequently
\[
T_f(\gamma_{\lambda_1}),\quad T_f(\gamma_{\lambda_2})\in\mathcal C_f
\]
have an incompatible TVG join, and $\mathcal C_f$ has no greatest member.
\end{theorem}
\begin{proof}
The map $\pi$ is continuous from each $(\Q,\gamma_{\lambda_i})$ to $(A,\lambda_i)$. It is therefore continuous between the joins. Thus $\psi\circ\pi$ is continuous on the rational join. If it were $\chi_s$ for $s\in\Q$, evaluation at $1$ would give $0=p(s)$ and hence $s\in\Z$. Surjectivity of $\pi$ would then give $\psi=s\iota$, a contradiction. Apply Lemma~\ref{lem:quotient-lift-v10}, Theorem~\ref{thm:transfer-v8}, and Proposition~\ref{prop:directed-v8}.
\end{proof}
No such pair $\lambda_1,\lambda_2$ has been constructed here. The implication is sufficient; its converse is not established. The torsion group $A$ is not a single Prüfer group: it is the direct sum, over primes $\ell$, of the groups $\Z(\ell^\infty)$, where $\Z(\ell^\infty)$ consists of the elements of $\Q/\Z$ whose orders are powers of $\ell$. The algebraic character group of $A$ is the profinite integer group $\widehat{\Z}$, the inverse limit of the finite cyclic groups $\Z/n\Z$. The forbidden extra characters above correspond to $\widehat{\Z}\setminus\Z$.

Arbitrary rational multiplication is not well-defined on $A$: the representatives $0$ and $1$ yield different classes after multiplication by $\frac12$. No such scalar structure is required on the quotient. Rational continuity is supplied by the full family defining $\gamma_\lambda$, and real scalar continuity is supplied by $T_f$.

Proposition~\ref{prop:rationalize-v8} equates existence of a greatest member of $\mathcal Q_\Gamma$ with existence of a Mackey topology for $(\Q,\Gamma)$. It does not identify $\mathcal Q_\Gamma$ with all of $\mathcal L_\Gamma$. On $\R$, Corollary~\ref{cor:necessary-v8} remains one-way. Failure of the unrestricted Mackey problem for an arbitrary character group must not be substituted for failure in the smaller TVG class $\mathcal C_f$.

\section{Polar criteria and the nuclearity restriction}\label{sec:polar-v10}
We give the precise duality conditions used from \cite{Gabriyelyan2022}. They provide tests for compatibility, rather than a resolution of the selected prime pair.

Let $G$ be an abelian group and $E\subseteq G^*$ a point-separating subgroup, where $G^*=\operatorname{Hom}(G,\T)$ contains all algebraic characters. Put $\T_+=[-\frac14,\frac14]+\Z$. For $U\subseteq G$ and $C\subseteq E$ define
\[
\begin{aligned}
\operatorname{Pol}_E(U)&=\{\chi\in E:\chi(U)\subseteq\T_+\},\\
\operatorname{Pol}_*(U)&=\{\chi\in G^*:\chi(U)\subseteq\T_+\},\\
\operatorname{Inv}(C)&=\{x\in G:\chi(x)\in\T_+\text{ for every }\chi\in C\}.
\end{aligned}
\]
These distinguish a polar inside the prescribed dual $E$ from one inside all algebraic characters. A subset $C\subseteq E$ is quasi-convex for the dual pair when
\[
C=\operatorname{Pol}_E(\operatorname{Inv}(C)).
\]
Write $\sigma(E,G)$ for pointwise convergence on $G$, and $\sigma(G,E)$ for the initial character topology on $G$. For $n\ge1$, $(n)C$ denotes the set of sums of $n$ elements of $C$, not merely the set of multiples $n\chi$.

\begin{definition}[The admissible compact sets]\label{def:mac-v10}
In the notation of \cite[Definition~3.3]{Gabriyelyan2022}, $\operatorname{Mac}(G,E)$ consists of the $\sigma(E,G)$-compact quasi-convex sets $C\subseteq E$ such that
\[
\operatorname{Pol}_*\bigl(\operatorname{Inv}((n)C+F)\bigr)\subseteq E
\]
for every $n\ge1$ and every finite subset $F\subseteq E$ containing zero.
\end{definition}
The extra polar condition is essential: compactness and quasi-convexity alone do not define this family. For a quasi-convex $C$, the polar topology $T_C$ has the zero-neighborhood base
\[
\operatorname{Inv}((n)C+F),\qquad n\ge1,\quad
F\subseteq E\text{ finite},\quad 0\in F.
\]
This is the polar group topology of \cite[Lemma~3.1]{Gabriyelyan2022}; it contains $\sigma(G,E)$.

\begin{theorem}[Established criteria specialized to the dual pair]\label{thm:mac-v10}
The topology $T_C$ is compatible with $E$ exactly when $C\in\operatorname{Mac}(G,E)$. Every compatible locally quasi-convex topology is a join of a nonempty directed family of these polar topologies. A greatest compatible locally quasi-convex topology exists exactly when $\operatorname{Mac}(G,E)$ is directed by inclusion.

For $C_1,C_2\in\operatorname{Mac}(G,E)$, compatibility of $T_{C_1}\vee T_{C_2}$ is equivalent to existence of some $C\in\operatorname{Mac}(G,E)$ containing $C_1\cup C_2$, and also to dual embeddedness of the diagonal in $(G,T_{C_1})\times(G,T_{C_2})$.
\end{theorem}
\begin{proof}
These are respectively Theorems~3.5, 3.8, 3.10 and Proposition~3.9 of \cite{Gabriyelyan2022}, applied to the dual pair $(G,E)$. The first condition also shows why the full algebraic polar is needed: a character bounded in the quarter arc on a basic neighborhood is continuous, so the condition excludes precisely additional characters.
\end{proof}
A subgroup is dually embedded when every continuous character on it extends to a continuous character of the ambient group. For our rational problem use $(G,E)=(\Q,\Gamma)$; for the quotient problem use $(A,\Z\iota)$. The theorem applies to these pairs, but its hypotheses do not automatically verify directedness. The representation of a topology may require an entire directed family, rather than one compact set. For the prime pair one must test common admissible upper bounds for the sets in the two representing families; a single unspecified compact set cannot replace that verification.

We use nuclear groups in the standard topological-group sense cited in \cite[Proposition~3.19]{Gabriyelyan2022}. The needed established facts are that locally compact abelian groups are nuclear, nuclear groups are closed under products and subgroups, and every subgroup of a nuclear group is dually embedded. These facts do not assert that every locally quasi-convex group is nuclear.

\begin{proposition}[Two nuclear members cannot be an incompatible pair]\label{prop:nuclear-v10}
Let $\delta_1,\delta_2$ be Hausdorff locally quasi-convex topologies on the same abelian group $G$, with common full character group $E$. If both $(G,\delta_1)$ and $(G,\delta_2)$ are nuclear, then $\delta_1\vee\delta_2$ is compatible with $E$.
\end{proposition}
\begin{proof}
The product is nuclear, so its diagonal is dually embedded. A continuous character of the diagonal extends to the product as $(x,y)\mapsto\chi_1(x)+\chi_2(y)$, where $\chi_i$ is $\delta_i$-continuous. Its restriction is $\chi_1+\chi_2\in E$. The diagonal topology is the join, and all members of $E$ remain continuous. Thus its full dual is $E$.
\end{proof}
Consequently an incompatible pair must have at least one non-nuclear member. This argument does not show that one nuclear member suffices for compatibility, nor that both members of a counterexample must be non-nuclear. A linear group topology means one with a zero-neighborhood base of subgroups. Such abelian topologies are nuclear by \cite[Proposition~3.19 and the following paragraph]{Gabriyelyan2022}; thus two linear members are excluded as well. Nuclearity of either of our prime topologies has not been established here. An embedding into a group of circle-valued null sequences with the supremum metric, by itself, does not settle nuclearity.

Theorem~3.20 of \cite{Gabriyelyan2022} ensures existence of a Mackey topology for locally quasi-convex abelian groups of finite exponent, meaning that one positive integer kills every element. Neither $\Q$ nor $A=\Q/\Z$ has finite exponent. That theorem therefore does not decide our cases. Also $\operatorname{Hom}(A,\T)$ is torsion-free: if $m\psi=0$, divisibility gives $x=my$ and $\psi(x)=m\psi(y)=0$. A finite-order extra-character mechanism cannot produce a nonzero character on $A$, just as Lemma~\ref{lem:roots-v8} excludes it on $\Q$ and $\R$.

The remaining tasks are to determine the dual of $\gamma_2\vee\gamma_3$, to check nuclearity rather than assume it, or to construct the quotient pair of Theorem~\ref{thm:quotient-pair-v10}. Failure of nuclearity alone would not prove incompatibility. Even a positive answer for the selected prime pair would not establish directedness for all compatible TVG topologies. These are limitations of the present arguments, not assertions that the wider literature has no answer.

\section{Real dilations and the Mackey problem for the Hamel pair}\label{sec:dilations-v13}
This section uses the real scalar structure of the class $\mathcal C_f$. Let $\mathcal L_f$ denote the class of all Hausdorff locally quasi-convex group topologies on the additive group $\R$ whose full continuous character group is $H_f$, with no assumption on scalar multiplication. Then $\mathcal C_f\subseteq\mathcal L_f$, and a greatest member of $\mathcal L_f$, if it exists, is a Mackey topology for the dual pair $(\R,H_f)$. For real $r\ne0$ write $m_r(x)=rx$. For a group topology $\eta$ on $\R$ put $\rho_r(\eta)=\{m_r^{-1}(U):U\in\eta\}$, so that $m_r:(\R,\rho_r(\eta))\to(\R,\eta)$ is an isomorphism of topological groups. Then $\rho_1(\eta)=\eta$ and $\rho_r(\rho_s(\eta))=\rho_{rs}(\eta)$. The \emph{continuous scalar ring} of $\eta$ is $C_\eta=\{r\in\R:m_r:(\R,\eta)\to(\R,\eta)\text{ is continuous}\}$. It is a subring of $\R$ containing $\Z$, because integer multiplication is continuous in every topological group, and $C_\eta=\R$ exactly when $\eta$ is a TVG topology.

\begin{lemma}[Automorphisms of the pair $(\R,H_f)$]\label{lem:aut-v13}
For every real $r\ne0$ one has $\{\chi\circ m_r:\chi\in H_f\}=H_f$. Conversely, if $\varphi:\R\to\R$ is an additive bijection such that $\chi\circ\varphi\in H_f$ for every $\chi\in H_f$, then $\varphi=m_d$ for some real $d\ne0$.
\end{lemma}
\begin{proof}
Since $\chi_{c,a}\circ m_r=\chi_{rc,ra}$ and $(c,a)\mapsto(rc,ra)$ is a bijection of the parameters, the first assertion holds. For the converse, fix a real number $s$. The ordinary character $x\mapsto p(sx)$ belongs to $H_f$, so $p(s\varphi(x))=p(c_sx+f(a_sx))$ for some reals $c_s,a_s$ and all $x$. The maps $x\mapsto s\varphi(x)$ and $x\mapsto c_sx+f(a_sx)$ are additive and real-valued with the same image in $\T$, so they coincide by Lemma~\ref{lem:lift}. Thus $s\varphi=c_s\,\mathrm{id}+f_{a_s}$. Taking $s=1$ and an irrational $s$ gives $\varphi=d\,\mathrm{id}+f_a$ and $s\varphi=e\,\mathrm{id}+f_b$, whence $(sd-e)x=f_b(x)-sf_a(x)$ for all $x$. The right side takes values in the countable set $\Q+s\Q$, whereas a nonzero multiple of the identity has uncountable image. Hence $sd=e$ and $sf_a=f_b$. Since $s$ is irrational and $f_a,f_b$ are rational-valued, $f_a=0$. Thus $\varphi=d\,\mathrm{id}$, and $d\ne0$ because $\varphi$ is bijective.
\end{proof}
\begin{proposition}[A greatest compatible topology is automatically a TVG topology]\label{prop:automatic-tvg-v13}
Suppose that $\mathcal L_f$ has a greatest member $\mu$. Then $\mu\in\mathcal C_f$, and $\mu$ is the greatest member of $\mathcal C_f$. Consequently the Greatest Compatible TVG Conjecture~\ref{conj:greatest-v12} holds for $f$ whenever the Mackey problem for the dual pair $(\R,H_f)$ has a positive answer.
\end{proposition}
\begin{proof}
Let $r\ne0$. The topology $\rho_r(\mu)$ is Hausdorff and locally quasi-convex, being isomorphic to $\mu$. A character $\chi$ is $\rho_r(\mu)$-continuous exactly when $\chi\circ m_r^{-1}$ is $\mu$-continuous, that is, when $\chi\in\{\psi\circ m_r:\psi\in H_f\}$, which equals $H_f$ by Lemma~\ref{lem:aut-v13}. Hence $\rho_r(\mu)\in\mathcal L_f$, so $\rho_r(\mu)\subseteq\mu$. This says that the identity $(\R,\mu)\to(\R,\rho_r(\mu))$ is continuous. Composing it with the isomorphism $m_r:(\R,\rho_r(\mu))\to(\R,\mu)$ shows that $m_r$ is continuous on $(\R,\mu)$. Since $r\ne0$ was arbitrary and $m_0$ is constant, $\mu$ is a TVG topology. Thus $\mu\in\mathcal C_f$, and because $\mathcal C_f\subseteq\mathcal L_f$ it is the greatest member of $\mathcal C_f$.
\end{proof}
The converse is not claimed: a greatest member of $\mathcal C_f$ is known to dominate only the compatible topologies that are themselves TVG topologies.

\begin{proposition}[Order of the dilated copies]\label{prop:order-v13}
Let $\eta$ be a group topology on $\R$ and let $r,t\ne0$. Then $\rho_r(\eta)\subseteq\rho_t(\eta)$ if and only if $r/t\in C_\eta$. Consequently the family $\{\rho_r(\eta):r\ne0\}$ is directed by inclusion if and only if the fraction field of the domain $C_\eta$ is $\R$.
\end{proposition}
\begin{proof}
The inclusion $\rho_r(\eta)\subseteq\rho_t(\eta)$ means that the identity $(\R,\rho_t(\eta))\to(\R,\rho_r(\eta))$ is continuous. Composing with the isomorphism $m_r:(\R,\rho_r(\eta))\to(\R,\eta)$, this holds exactly when $m_r:(\R,\rho_t(\eta))\to(\R,\eta)$ is continuous. Since $m_r=m_{r/t}\circ m_t$ and $m_t:(\R,\rho_t(\eta))\to(\R,\eta)$ is an isomorphism, this is equivalent to continuity of $m_{r/t}$ on $(\R,\eta)$, that is, to $r/t\in C_\eta$.

Suppose the family is directed. For $r\ne0$ choose $t$ with $\rho_r(\eta)\subseteq\rho_t(\eta)$ and $\rho_1(\eta)\subseteq\rho_t(\eta)$. Then $r/t,1/t\in C_\eta$ and $r=(r/t)/(1/t)$ lies in the fraction field. Conversely, suppose every real number is a quotient of elements of $C_\eta$. Given $r,s\ne0$ write $r=a/b$ and $s=c/d$ with $a,b,c,d\in C_\eta\setminus\{0\}$ and put $t=1/(bd)$. Then $r/t=ad$ and $s/t=cb$ lie in $C_\eta$, so $\rho_r(\eta)$ and $\rho_s(\eta)$ are both contained in $\rho_t(\eta)$.
\end{proof}
\begin{corollary}[Real saturation]\label{cor:real-saturation-v13}
Let $\eta\in\mathcal L_f$ and suppose that the fraction field of $C_\eta$ is $\R$. Then $\eta^{[\R]}=\bigvee_{r\ne0}\rho_r(\eta)$ belongs to $\mathcal C_f$, and it is the least TVG topology containing $\eta$.
\end{corollary}
\begin{proof}
As in the proof of Proposition~\ref{prop:automatic-tvg-v13}, each $\rho_r(\eta)$ lies in $\mathcal L_f$. By Proposition~\ref{prop:order-v13} the family is directed. The join of a nonempty directed family in $\mathcal L_f$ lies in $\mathcal L_f$, by the proof of Proposition~\ref{prop:directed-v8}, which uses only the quarter-arc criterion and directedness. Let $s,t\ne0$. The map $m_s:(\R,\rho_{ts}(\eta))\to(\R,\rho_t(\eta))$ is continuous, because $m_t\circ m_s=m_{ts}$ is an isomorphism of $(\R,\rho_{ts}(\eta))$ onto $(\R,\eta)$ and $m_t$ is an isomorphism of $(\R,\rho_t(\eta))$ onto $(\R,\eta)$. Since $\rho_{ts}(\eta)\subseteq\eta^{[\R]}$, the map $m_s$ is continuous from $(\R,\eta^{[\R]})$ into every $\rho_t(\eta)$, hence into their join. So $\eta^{[\R]}$ is a TVG topology and lies in $\mathcal C_f$. Finally let $\nu\supseteq\eta$ be a TVG topology. For $r\ne0$ the identity $(\R,\nu)\to(\R,\rho_r(\eta))$ is continuous exactly when $m_r:(\R,\nu)\to(\R,\eta)$ is; this holds because $m_r$ is continuous on $(\R,\nu)$ and $\nu\supseteq\eta$. Hence $\nu\supseteq\rho_r(\eta)$ for all $r$, and $\nu\supseteq\eta^{[\R]}$.
\end{proof}
For rational dilations the copies pulled back along $x\mapsto x/n!$ form an increasing cofinal chain among the nonzero rational dilations: integer multiplication is always continuous. Thus no additional scalar-ring condition is needed for rational saturation. For real dilations the criterion of Proposition~\ref{prop:order-v13} is exact. If $C_\eta$ is countable, for instance $C_\eta=\Q$, its fraction field is countable and different from $\R$, the family is not directed, and the directed-join argument gives no conclusion about compatibility of the join over the real orbit. Failure of directedness does not itself imply incompatibility.

\begin{proposition}[A necessary condition for a greatest member]\label{prop:necessary-v13}
Suppose that $\mathcal C_f$ has a greatest member $\mu$. Then $\tu\subseteq\mu$, and for every $\eta\in\mathcal C_f$ one has $\eta\vee\tu=\eta\vee\sigma_f$, with full character group $H_f$. Consequently a TVG topology $\eta\in\mathcal C_f$ for which $\eta\vee\sigma_f$ has a continuous character outside $H_f$ would show that $\mathcal C_f$ has no greatest member. By Proposition~\ref{prop:nuclear-v10} such an $\eta$ would not be nuclear, since $\sigma_f$ is nuclear.
\end{proposition}
\begin{proof}
We have $\sigma_f\in\mathcal C_f$ and $\tu\subseteq\sigma_f$, so $\tu\subseteq\sigma_f\subseteq\mu$. Every $\eta\in\mathcal C_f$ contains the weak topology $\beta_f$ generated by $H_f$. Since $\sigma_f=\tu\vee\beta_f$, the two joins are equal. Both are contained in $\mu$, so each of their continuous characters is $\mu$-continuous and belongs to $H_f$. Every member of $H_f$ is $\eta$-continuous, hence continuous for the joins. The last assertion uses the facts recalled before Proposition~\ref{prop:nuclear-v10}: $(\R,\sigma_f)$ is a subgroup of the locally compact group $\R_{\mathrm{usual}}\times K$ and is therefore nuclear, so if $\eta$ were also nuclear the join would be compatible.
\end{proof}
\begin{corollary}[Position of the conjecture]\label{cor:sandwich-v13}
For every nonzero $f\in D$, the following implications hold: a Mackey topology for $(\R,H_f)$ exists $\Longrightarrow$ Conjecture~\ref{conj:greatest-v12} holds for $f$ $\Longrightarrow$ a Mackey topology for $(\Q,\Gamma)$ exists.
\end{corollary}
\begin{proof}
The first implication is Proposition~\ref{prop:automatic-tvg-v13}. The second is Corollary~\ref{cor:necessary-v8}.
\end{proof}
Neither converse is proved here. A compatible real saturation for every member of $\mathcal L_f$ would suffice for the first converse, but is not established. For the second converse, one sufficient route would extend the rational Mackey topology to a cofinal family in $\mathcal C_f$; no such argument is provided. Neither statement asserts equivalence with the Mackey problem for $\Q/\Z$.

\section{Conclusion}\label{sec:answers}
The three Hamel constructions use different observations of the same coefficient function. The topology $\tau_f$ observes real-valued coefficients, $\beta_f$ observes ordinary and Hamel characters, and $\sigma_f$ adjoins the usual topology to $\beta_f$. Their definitions and neighborhood bases are given in Section~4. The investigation yields the following conclusions.

\begin{enumerate}
\item \textbf{Comparison and failure of real-valued representation.}
For every nonzero $f\in D$,
\[
\tu\subsetneq\sigma_f\subsetneq\tau_f.
\]
Every nonzero rational-valued additive function is unbounded on every $\sigma_f$-neighborhood of zero. Hence modular observations together with the usual topology do not recover real-valued coefficient continuity, and $\sigma_f$ cannot be generated by usual-rational-valued additive observations (Theorem~\ref{thm:main} and Corollary~\ref{cor:representation}).

\item \textbf{Characters and completions.}
The topologies $\beta_f$ and $\sigma_f$ have exactly the same characters $\chi_{c,a}(x)=p(cx+f(ax))$, with unique parameters. Their continuous real-valued additive functions are respectively zero and the scalar multiples of the identity. Their completions are $K\times K$ and $\R_{\mathrm{usual}}\times K$, where $K=\widehat{\R_d}$ (Theorems~\ref{thm:characters}, \ref{thm:real-duals}, and \ref{thm:completions}). Both original groups are totally disconnected, incomplete, nonmetrizable, and not locally compact, although their completions are connected.

\item \textbf{Rigidity and isomorphism types.}
Every continuous additive map between any two modular-family members is a scalar multiplication. Within either family, a nonzero such map exists precisely when the topologies coincide, equivalently when their coefficient functions differ by a nonzero argument rescaling. Each family therefore supplies $2^{\mathfrak c}$ pairwise nonisomorphic topological groups with the same completion model (Theorem~\ref{thm:rigidity} and Corollary~\ref{cor:isomorphism}). No classification by arbitrary homeomorphisms is asserted.

\item \textbf{Compactness and convergence.}
The topologies $\beta_f$ and $\sigma_f$ have the same compact subsets, the same induced topologies on those subsets, the same convergent sequences, and the same compact-open topology on their common character set (Theorem~\ref{thm:compact}). Explicit compactness and convergence tests also distinguish what must be checked for $\tau_f$ (Corollary~\ref{cor:sequences} and Proposition~\ref{prop:compact-tests}). The universal sequence $x_n=n!(3-2\sqrt2)^{n^2}$ converges in $\sigma_f$ and $\beta_f$ but not in $\tau_f$ for every nonzero $f$. Its terms together with zero form a compact set for the modular pair but not for $\tau_f$ (Theorem~\ref{thm:universal-sequence} and Corollary~\ref{cor:compact-separation}). Theorem~\ref{thm:finite-span-sequences} gives an exact classification in finite rational spans; dimension two is the least possible dimension for separation.

\item \textbf{Precompact classification and minimality.}
All Hausdorff precompact TVG topologies correspond exactly to nonzero real subspaces of the algebraic character space under argument rescaling. There are $2^{2^{\mathfrak c}}$ such topologies and $2^{\mathfrak c}$ minimal ones within this class. The minimal topology $t$ and the supplementary minimal topology $\alpha_f$ have no common coarser Hausdorff TVG topology. Thus no least Hausdorff TVG topology exists (Theorem~\ref{thm:precompact} and Corollary~\ref{cor:minimal}). Minimality of $t$ remains fully consistent with this conclusion.

\item \textbf{Compatibility and maximal nondiscrete extensions.}
For the common character set $H_f$, $\beta_f$ is least compatible, and $\sigma_f$ is least compatible among topologies containing the usual topology (Proposition~\ref{prop:compatible}). Using the topology on $\Q/\Z$ furnished by de la Barrera Mayoral \cite{Barrera2014}, Theorem~\ref{thm:finer-compatible} constructs a strictly finer compatible locally quasi-convex TVG refinement $\eta_f$ of $\sigma_f$. Thus $\sigma_f$ is not finest. Every nondiscrete Hausdorff TVG admits a maximal nondiscrete extension. Every maximal member contains $\tu$ and the full discrete-group Bohr topology $B$, makes every algebraic character continuous, and has only finite compact subsets and eventually constant convergent sequences (Theorem~\ref{thm:maximal}). There is no greatest nondiscrete TVG, even among topologies finer than usual: the topologies $\rho_E$ and $\rho_F$ are nondiscrete but have discrete join. Their maximal extensions are different, and both extend $\tu\vee B$ (Theorem~\ref{thm:no-greatest} and Corollary~\ref{cor:nonunique-maximal}).
\item \textbf{Rational transfer and ultrafilter criteria.}
Every member of $\mathcal C_f$ has maximal compatible extensions. A greatest member exists exactly when every binary join stays compatible (Proposition~\ref{prop:directed-v8}). Rational topologies with full dual $\Gamma$ transfer by the coefficient maps into $\mathcal C_f$, and incompatible rational joins would give incompatible TVG joins (Theorem~\ref{thm:transfer-v8}). A greatest member would imply existence of a Mackey topology for $(\Q,\Gamma)$; the converse is not established (Corollary~\ref{cor:necessary-v8}). The ultrafilter results give exact continuity tests and a sufficient counterexample test, without producing the needed neighborhood control.

\item \textbf{Explicit candidates and the unsolved step.}
The topologies generated by $P_{2,d}$ and $P_{3,d}$ have full character set $\Gamma$, are metrizable and incomparable, and are not locally precompact (Theorems~\ref{thm:prime-v8} and \ref{thm:prime-incomparable-v8}). Their join remains strictly weaker than usual (Proposition~\ref{prop:prime-join-v8}). Its compatibility is equivalent to a torsion condition on an explicit integer character group (Proposition~\ref{prop:integer-v8}); that condition is not proved. An alternating dyadic block construction gains irrational characters already in the individual topologies and therefore fails compatibility (Proposition~\ref{prop:blocks-v8}). Appendix~\ref{app:segments-v12} gives a second obstruction: the full factorial-segment topology on $\Q/\Z$ is discrete, and its raw odd-indexed half is not Hausdorff. These failed constructions and the known finite-order mechanism do not prove nonexistence for $\mathcal C_f$.

\item \textbf{Quotient lifting with the TVG property preserved.}
Every compatible locally quasi-convex topology $\lambda$ on $\Q/\Z$ lifts to $\gamma_\lambda\in\mathcal Q_\Gamma$ (Lemma~\ref{lem:quotient-lift-v10}). An incompatible quotient pair would yield an incompatible rational pair and hence two compatible TVG topologies on $\R$ with incompatible join (Theorem~\ref{thm:quotient-pair-v10}). This proves a sufficient reduction, not existence of the required pair or the converse implication. Arbitrary joins preserve scalar continuity, while their compatibility remains an independent condition (Proposition~\ref{prop:scalar-join-v10}).

\item \textbf{Polar criteria and the nuclearity obstruction.}
The full definition of $\operatorname{Mac}(G,E)$ includes an algebraic-polar condition as well as weak compactness and quasi-convexity. The established directedness and diagonal-extension criteria are recorded with their source in Theorem~\ref{thm:mac-v10}. If both members of a compatible pair are nuclear, their join is compatible (Proposition~\ref{prop:nuclear-v10}); at least one member of an incompatible pair must therefore be non-nuclear. This does not establish nuclearity of either prime topology or compatibility of their join.

\item \textbf{Real dilations and the dyadic dichotomy.}
The only additive automorphisms of the dual pair $(\R,H_f)$ are the real scalars (Lemma~\ref{lem:aut-v13}). Hence a greatest member of the unrestricted class $\mathcal L_f$, if it exists, is automatically a TVG topology and the greatest member of $\mathcal C_f$ (Proposition~\ref{prop:automatic-tvg-v13}); so the Greatest Compatible TVG Conjecture follows from the Mackey problem for $(\R,H_f)$ and itself implies a Mackey topology for $(\Q,\Gamma)$ (Corollary~\ref{cor:sandwich-v13}). The dilated copies of a compatible topology $\eta$ are directed exactly when the fraction field of its continuous scalar ring is $\R$, and then a real saturation exists (Proposition~\ref{prop:order-v13} and Corollary~\ref{cor:real-saturation-v13}). A greatest member of $\mathcal C_f$ forces $\eta\vee\tu$ and $\eta\vee\sigma_f$ to be compatible for every $\eta\in\mathcal C_f$ (Proposition~\ref{prop:necessary-v13}). No splitting of the dyadic chain gives an incompatible pair of the topologies $\xi_S$ (Theorem~\ref{thm:dyadic-v13}).

\end{enumerate}

The principal extremal question still unresolved here is whether the character set $H_f$ admits a greatest compatible locally quasi-convex TVG topology. Its proposed affirmative form is the Greatest Compatible TVG Conjecture~\ref{conj:greatest-v12}, explicitly marked unproved. Theorem~\ref{thm:finer-compatible} shows that $\sigma_f$ cannot be that member. Sections~\ref{sec:joins-v8}--\ref{sec:dilations-v13} record the proved reductions and the unresolved mixed-prime character calculation, while Appendix~\ref{app:blocks-v11} documents the unsuccessful alternating-block attempt; Proposition~\ref{prop:conditional-v8} is conditional, not a nonexistence theorem for this compatible class. The quotient lifting and the nuclearity obstruction add precise sufficient conditions but do not decide that extremal question. Sequence separation and nonexistence of a greatest nondiscrete TVG have been proved above. The sequence classification in finite rational spans does not claim a classification of every infinite-dimensional rational-span sequence. These scope statements concern what has been proved here and make no assertion that the wider literature lacks further answers.

The initial-topology mechanism, Bohr construction, compact duality, Glicksberg theorem, and compatible-topology framework are standard background. The results here specialize those tools and provide direct arguments for the stated families. Historical priority for specialized statements remains unverified, including comparison with Hejcman's full paper.

\appendix
\section{An unsuccessful alternating block construction}\label{app:blocks-v11}
This appendix documents a failed attempt to obtain compatible topologies with an incompatible join. Its obstruction is proved below: irrational characters are already continuous in each individual topology. The construction is included to explain why this route does not settle the questions in Section~\ref{sec:status-v8}.

Another concrete attempt partitions the nonnegative integers into alternating blocks
\[
B_k=\{k^2,\ldots,(k+1)^2-1\},\quad
I=\bigcup_{j\ge0}B_{2j},\quad
J=\bigcup_{j\ge0}B_{2j+1}.
\]
For $S=I,J$, let $\xi_S$ be generated on $\Q$ by
\[
Q_{S,d}(q)=\sup_{n\in S}
\left\|\frac{q}{d\,2^n}\right\|_{\T},
\qquad d\in\mathbb N_{>0}.
\]
These are metrizable locally quasi-convex topologies with rational scalar continuity. Their join is usual: $I\cup J$ contains all nonnegative integers, and for $|q|\ge\frac{1}{4}$ some $n$ satisfies $\frac{1}{4}\le \frac{|q|}{2^n}<\frac{1}{2}$. Thus the ball of radius $\frac{1}{4}$ for the supremum over all $n$ is $(-\frac{1}{4},\frac{1}{4})\cap\Q$.

\begin{proposition}[Why the block attempt fails]\label{prop:blocks-v8}
Both $\xi_I$ and $\xi_J$ already have irrational-parameter continuous characters and hence do not belong to $\mathcal Q_\Gamma$.
\end{proposition}
\begin{proof}
For a block $[a,b]\subseteq S$, suppose $Q_{S,1}(q)<\delta<\frac{1}{4}$. Choose the centered representative $y_n$ of $\frac{q}{2^n}+\Z$. For consecutive indices in the block, $y_n-2y_{n+1}$ is an integer of absolute value less than $3\delta<1$, so it is zero. Therefore
\[
\|\frac{q}{2^b}\|_{\T}\le\delta\,2^{-(b-a)}.
\]
The endpoints of the $I$-blocks and $J$-blocks give respectively
\[
\alpha_I=\sum_{j\ge0}2^{-(4j^2+4j)},\qquad
\alpha_J=\sum_{j\ge0}2^{-(4j^2+8j+3)}.
\]
Both numbers are irrational: their binary expansions have infinitely many isolated ones with unbounded gaps and therefore are not eventually periodic. Summing the centered representatives over block endpoints gives
\[
\|\alpha_Iq\|_{\T}\le\frac{16}{15}\delta,
\qquad
\|\alpha_Jq\|_{\T}\le\frac{4}{15}\delta
\]
on the corresponding balls. The circle series agrees with the real series modulo $\Z$, and the displayed absolute bounds justify its convergence. These estimates prove continuity of the two irrational characters in their respective topologies.
\end{proof}
The block attempt obtains a join with extra characters but fails individual compatibility. The prime-square construction obtains individual compatibility and incomparability but has no computed join dual. These two outcomes must not be conflated.

\section{Factorial segments and finite-search limitations}\label{app:segments-v12}
We consider a second unsuccessful family on $A=\Q/\Z$, with inclusion $\iota:A\to\T$ as in Section~\ref{sec:scalar-v10}. All scalar coefficients in this appendix are integers acting on the torsion group. No multiplication of a class in $A$ by an arbitrary rational number is being defined.

For a set $C\subseteq\Z$ containing $0$ and $1$, put
\[
P_C(x)=\sup_{c\in C}\|c x\|_{\T},\qquad
U_C(\varepsilon)=\{x\in A:P_C(x)<\varepsilon\}.
\]
Let $\lambda_C$ be the topology of uniform convergence on these characters, generated by these balls. It is a Hausdorff locally quasi-convex group topology: the observation $1$ separates points, the circle triangle inequality gives the group property, and the closed balls of radii $\frac{1}{4n}$ are intersections of inverse images of quasi-convex arcs. It contains the circle topology. We call $C$ compatible for this construction when the full character group of $(A,\lambda_C)$ is $\Z\iota$. Compactness of $C$ alone does not imply that compatibility.

The profinite topology on $\Z$ has the subgroups $m\Z$, $m\ge1$, as a zero-neighborhood base. Here a set is profinitely null when, for each $m\ge1$, all but finitely many of its elements lie in $m\Z$. Adjoining zero to such a set makes it compact: any open cover has a member containing zero and all but finitely many points, and the remaining finitely many points have a finite subcover. Adding a finite set, such as $\{1\}$, preserves this property.

Choose finite integers $J_k\ge k$ for $k\ge1$, and define
\[
C=\{0\}\cup\{k!\,2^j:k\ge1,\ 0\le j\le J_k\}.
\]
This set contains $1$. It is profinitely null, since $m$ divides $k!$ for all $k\ge m$, and the finitely many earlier segments are finite. The following estimate explains its topology.

\begin{lemma}[Contraction along a finite doubling segment]\label{lem:segment-v12}
Let $t\in\T$, let $J\ge0$ be an integer, and suppose
\[
\|2^jt\|_{\T}<\varepsilon
\qquad(0\le j\le J),\qquad 0<\varepsilon<\frac18.
\]
Then $\|t\|_{\T}<\varepsilon\,2^{-J}$.
\end{lemma}
\begin{proof}
Choose centered real representatives $u_j$ of $2^jt$ with $|u_j|<\varepsilon$. The number $2u_j-u_{j+1}$ is an integer of absolute value less than $3\varepsilon<1$, so it is zero. Thus $u_J=2^Ju_0$, which gives the estimate.
\end{proof}

\begin{proposition}[The factorial-segment topology is discrete]\label{prop:segments-discrete-v12}
For the set $C$ above and every $0<\varepsilon<\frac18$,
\[
U_C(\varepsilon)=\{0\}.
\]
Consequently $\lambda_C$ is discrete and is not compatible with $\Z\iota$.
\end{proposition}
\begin{proof}
Suppose $x\ne0$ belongs to $U_C(\varepsilon)$. Since $x$ has finite order, there is a least integer $m\ge2$ with $m!x=0$. Put $t=(m-1)!x$. Then $t\ne0$ and $mt=0$, so, under the inclusion into $\T$,
\[
\|t\|_{\T}\ge\frac1m.
\]
Apply Lemma~\ref{lem:segment-v12} to the segment with $k=m-1$. Since $J_{m-1}\ge m-1$ and $m\le2^{m-1}$, it gives
\[
\|t\|_{\T}
<\varepsilon\,2^{-J_{m-1}}
\le\varepsilon\,2^{-(m-1)}
\le\frac{\varepsilon}{m}
<\frac1m,
\]
a contradiction. Thus the indicated ball contains only zero. The discrete group $A$ has all its algebraic characters continuous, with character group $\widehat{\Z}$, strictly larger than $\Z\iota$ as described in Section~\ref{sec:scalar-v10}.
\end{proof}
This proves an obstruction for every choice of finite $J_k\ge k$. For the full infinite segment family there are no nonzero points hidden beyond a finite search bound in these small balls. Truncating the family changes the topology; a numerical computation with capped segments cannot establish the full assertion by itself.

\begin{proposition}[The raw odd-indexed half is not Hausdorff]\label{prop:odd-segments-v12}
Let $C_{\mathrm{odd}}$ consist of the observations $k!\,2^j$ with $j$ odd. The topology generated by these observations alone is not Hausdorff. In particular, the character
\[
\psi(x)=\sum_{k\ge1}k!x
\]
is not continuous for that raw topology.
\end{proposition}
\begin{proof}
Every observation is an even integer and vanishes at $x=\frac12+\Z$. This nonzero point is therefore topologically indistinguishable from zero. The displayed sum defines an algebraic character: for any fixed $x\in A$ only finitely many terms are nonzero, and for any two arguments a common finite cutoff verifies additivity. At $\frac12+\Z$ its value is $\frac12+\Z$, since all terms with $k\ge2$ vanish. A continuous character into the Hausdorff circle must vanish on the closure of zero, proving the assertion.
\end{proof}
For a Hausdorff candidate one must include the circle topology, for example by adjoining the observation $1$. That change does not by itself prove compatibility. The same profinite-series character cannot be asserted continuous in the raw odd-indexed half merely by repeating the estimate for the full family.

\textbf{The quantifiers in a numerical continuity test.}
For any algebraic character $\chi:A\to\T$, the quarter-arc criterion gives
\[
\begin{aligned}
\chi\text{ is }\lambda_C\text{-continuous}
\quad&\Longleftrightarrow\quad
\text{there exists }\varepsilon>0\text{ such that }
\|\chi(x)\|_{\T}\le\frac14\text{ for all }x\in U_C(\varepsilon).
\end{aligned}
\]
Discontinuity instead requires a violating point in every such neighborhood, or in a sequence of these neighborhoods whose radii decrease to zero. A finite list of violations at selected radii does not prove discontinuity; passing a bounded search does not prove continuity. No empirical output is used as a proof in this appendix.

\begin{question}[Unions of profinitely null compact sets]\label{question:null-unions-v12}
Is the following proposed closure assertion true? If $C_1,C_2\subseteq\Z$ are profinitely null compact sets containing $0$ and $1$, and each $\lambda_{C_i}$ has full character group $\Z\iota$, then $\lambda_{C_1\cup C_2}$ also has full character group $\Z\iota$.
\end{question}
This statement is unresolved here. The identity
\[
P_{C_1\cup C_2}=\max\{P_{C_1},P_{C_2}\}
\]
gives $\lambda_{C_1\cup C_2}=\lambda_{C_1}\vee\lambda_{C_2}$. Thus the question asks about compatibility of joins within this restricted construction. Its contrapositive says that every such decomposition of an incompatible union has at least one incompatible member; the quantifier is over every decomposition, not existence of one bad part.

A positive answer would exclude counterexamples of this restricted kind. To infer existence of a Mackey topology for $(A,\Z\iota)$, a further argument would have to show that this construction is cofinal among all compatible locally quasi-convex topologies, or otherwise establish the full directedness condition in Theorem~\ref{thm:mac-v10}. No such reduction is proved here. A negative answer, together with the two individual compatibility proofs, would give the quotient pair required by Theorem~\ref{thm:quotient-pair-v10}. Neither the discrete full segment family nor its unverified splits supplies that pair.

\section{Dyadic subset splittings cannot give an incompatible pair}\label{app:dyadic-v13}
Appendix~\ref{app:blocks-v11} treated one splitting of the dyadic chain $\{2^n\}$. We show that no splitting of this type works. For $S\subseteq\mathbb N_0$ and an integer $d\ge1$ put $Q_{S,d}(q)=\sup_{n\in S}\|q/(d\,2^n)\|_{\T}$ for $q\in\Q$, with the supremum over the empty set equal to $0$. Let $\xi_S$ be the group topology on $\Q$ with zero-neighborhood base $\{q:Q_{S,d}(q)<\varepsilon\}$, $d\ge1$, $\varepsilon>0$. A \emph{run} in $S$ is a set of consecutive integers contained in $S$. If $S$ is infinite, its \emph{gaps} are the differences of consecutive elements.

\begin{theorem}[Dichotomy for dyadic subsets]\label{thm:dyadic-v13}
Let $S\subseteq\mathbb N_0$.

(1) If $S$ is infinite and has bounded gaps, then $\xi_S$ is the usual topology of $\Q$. Its character group $\{q\mapsto p(cq):c\in\R\}$ is strictly larger than $\Gamma$.

(2) If $S$ contains runs of arbitrary length, then $\xi_S$ has a continuous character $q\mapsto p(\gamma q)$ with $\gamma$ irrational. In particular its character group differs from $\Gamma$.

(3) For every $S$, at least one of $\xi_S$ and $\xi_{\mathbb N_0\setminus S}$ has character group different from $\Gamma$, and therefore does not belong to $\mathcal Q_\Gamma$.
\end{theorem}
\begin{proof}
(1) Let $n_0=\min S$, let $L$ bound the gaps, and let $0<\varepsilon<2^{-L-2}$. Suppose $Q_{S,1}(q)<\varepsilon$ and let $y_n$ be the centered representative of $q/2^n$ modulo $\Z$, for $n\in S$. For consecutive $n<n'$ in $S$ with $g=n'-n\le L$ we have $y_n\equiv 2^gy_{n'}$ modulo $\Z$, while $|y_n|<\varepsilon$ and $2^g|y_{n'}|<2^L\varepsilon<\frac14$. The integer $y_n-2^gy_{n'}$ therefore has absolute value less than $\frac12$, so it is zero and $y_n=2^gy_{n'}$. For $n\in S$ large enough that $|q|/2^n<\frac12$ we have $y_n=q/2^n$, and descending along consecutive elements of $S$ gives $y_n=q/2^n$ for every $n\in S$. In particular $|q|<\varepsilon\,2^{n_0}$. Conversely, if $|q|<\varepsilon\,2^{n_0}$ then $|q|/2^n<\varepsilon<\frac12$ for every $n\in S$, so $Q_{S,1}(q)<\varepsilon$. Thus $\{q:Q_{S,1}(q)<\varepsilon\}=(-\varepsilon2^{n_0},\varepsilon2^{n_0})\cap\Q$. Since $Q_{S,d}(q)=Q_{S,1}(q/d)$, the neighborhood base of $\xi_S$ is a base of the usual topology of $\Q$. Continuous characters of the usual rationals extend to $\R$, so they are the maps $q\mapsto p(cq)$ with $c\in\R$, and these exceed $\Gamma$.

(2) Choose runs $[a_i,b_i]\subseteq S$ for $i\ge1$ with $b_i-a_i\ge i+1$ and $a_{i+1}\ge b_i+i$. Put $\gamma=\sum_{i\ge1}2^{-b_i}$. The binary expansion of $\gamma$ has ones exactly at the positions $b_i$, and the gaps $b_{i+1}-b_i\ge i$ are unbounded, so the expansion is not eventually periodic and $\gamma$ is irrational. Let $\delta=\frac15$ and suppose $Q_{S,1}(q)<\delta$. For $n\in S$ let $y_n$ be the centered representative of $q/2^n$. For consecutive indices $n,n+1$ in a run, $y_n-2y_{n+1}$ is an integer of absolute value less than $3\delta<1$, hence zero. Therefore $|y_{b_i}|=2^{-(b_i-a_i)}|y_{a_i}|<\delta\,2^{-(i+1)}$. Because $q2^{-b_i}\to0$ and $y_{b_i}\to0$, the integers $q2^{-b_i}-y_{b_i}$ vanish for large $i$, so $\gamma q\equiv\sum_iy_{b_i}$ modulo $\Z$. Hence $\|\gamma q\|_{\T}<\delta\sum_i2^{-(i+1)}=\delta/2<\frac1{10}$. The neighborhood $\{q:Q_{S,1}(q)<\frac15\}$ is mapped into the arc $[-\frac1{10},\frac1{10}]+\Z\subseteq[-\frac14,\frac14]+\Z$, so $q\mapsto p(\gamma q)$ is continuous by the quarter-arc criterion. It does not belong to $\Gamma$: if $p(\gamma q)=p(rq)$ for all $q\in\Q$ with $r\in\Q$, then $(\gamma-r)q\in\Z$ for all $q$, which forces $\gamma=r$.

(3) If $S$ is infinite with bounded gaps, apply (1) to $S$. If $S$ is infinite with unbounded gaps, then $\mathbb N_0\setminus S$ contains runs of arbitrary length, and (2) applies to $\xi_{\mathbb N_0\setminus S}$. If $S$ is finite, then $\mathbb N_0\setminus S$ is infinite with bounded gaps, and (1) applies to it.
\end{proof}
This theorem excludes every splitting $\{2^n:n\in S\}\cup\{2^n:n\notin S\}$ of the dyadic chain. It does not exclude other chains, nor non-chain compact sets of denominators. For finite nonempty $S$ the topology $\xi_S$ is the precompact topology generated by the characters $q\mapsto q/(d2^n)+\Z$, so it is compatible; its partner $\xi_{\mathbb N_0\setminus S}$ is the usual topology.

\section*{Acknowledgement}
\addcontentsline{toc}{section}{Acknowledgement}
The author acknowledges OpenAI's Codex and Anthropic's Claude for assistance in mathematical discussion, review of arguments, exposition, and manuscript preparation. The author assumes responsibility for the final content and any remaining errors.


\addcontentsline{toc}{section}{References}
\begin{thebibliography}{99}
\raggedright
\bibitem{Hamel1905} G. Hamel, Eine Basis aller Zahlen und die unstetigen L\"osungen der Funktionalgleichung $f(x+y)=f(x)+f(y)$, Mathematische Annalen 60 (1905), 459--462.
\bibitem{Y1975} I. B. Yaroshevskaya, An example of a quasi-vector topology on $\R^1$ other than the usual topology, Vestnik Moskov. Univ. Ser. I Mat. Mekh., no. 2 (1975), 35--40. English translation: An example of a separable quasi-vector topology in $\R^1$ different from the usual topology, Moscow University Mathematics Bulletin 30, nos. 1--2 (1975), 100--104. Zbl 0307.54038.
\bibitem{Y1981} I. B. Yaroshevskaya, On the minimal topology in the class of topological vector groups on $\R^1$, Russian Mathematical Surveys 36:6 (1981), 205--206; Russian original: Uspekhi Matematicheskikh Nauk 36:6(222) (1981), 235--236. \href{https://doi.org/10.1070/RM1981v036n06ABEH003139}{doi:10.1070/RM1981v036n06ABEH003139}.
\bibitem{D1978} V. M. Dergachev, Topological vector groups, Russian Mathematical Surveys 33:4 (1978), 251--252; Russian original: Uspekhi Matematicheskikh Nauk 33:4(202) (1978), 209--210.
\bibitem{DT2008} X. Dom\'inguez and V. Tarieladze, Group topologies on vector spaces and character lifting properties, Bolet\'in de la Sociedad Matem\'atica Mexicana (3) 14:1 (2008), 21--34. \href{https://caminos.udc.es/gmni/pdf/2008/group_topologies_on_vector_spaces_and_character_lifting_property.pdf}{Online copy}.
\bibitem{H1988} J. Hejcman, Topological vector group topologies for the real line, in Z. Frol\'ik (ed.), General Topology and Its Relations to Modern Analysis and Algebra VI, Research and Exposition in Mathematics 16, Heldermann, Berlin, 1988, 241--248. MR 0952610.
\bibitem{Y2026} I. Yaroshevskaya, Topological Vector Group Topologies Between the Minimal Topology and the Usual Topology on the Real Line, arXiv:2609.26668 [math.GN] (2026). \href{https://arxiv.org/abs/2609.26668}{arXiv:2609.26668}.
\bibitem{YHamel} I. Yaroshevskaya, Hamel-based TVG Stronger than the Usual Topology on $\R$, arXiv:2610.24132 (2026). \href{https://arxiv.org/abs/2610.24132}{arXiv:2610.24132}.

\bibitem{Morris1977} S. A. Morris, Pontryagin Duality and the Structure of Locally Compact Abelian Groups, London Mathematical Society Lecture Note Series 29, Cambridge University Press, Cambridge, 1977. \href{https://www.cambridge.org/core/books/pontryagin-duality-and-the-structure-of-locally-compact-abelian-groups/854007D98ACCC770CCD9432D3D5D68D0}{Publisher record}.
\bibitem{Elkies2002} N. D. Elkies, Outline of Pontrjagin duality, Math 55b course notes, Harvard University (2002). \href{https://people.math.harvard.edu/~elkies/M55b.02/pontrjagin.html}{Online notes}.
\bibitem{EW2010} M. Einsiedler and T. Ward, Ergodic Theory with a View towards Number Theory, draft dated 21 June 2010, Appendix C: Topological Groups. \href{https://webspace.maths.qmul.ac.uk/f.vivaldi/teaching/ETAD/NotesI.pdf}{Online draft}.
\bibitem{MegNotes} M. Megrelishvili, Lecture Notes in Topological Groups (2014), online revision dated 7 February 2017, Theorem 3.7, Remark 6.11, and Theorem 8.3. \href{https://u.math.biu.ac.il/~megereli/TGrNotes070217.pdf}{Online notes}.

\bibitem{Glicksberg1962} I. Glicksberg, Uniform boundedness for groups, Canadian Journal of Mathematics 14 (1962), 269--276. \href{https://doi.org/10.4153/CJM-1962-017-3}{doi:10.4153/CJM-1962-017-3}.
\bibitem{Gabriyelyan2022} S. Gabriyelyan, Compatible group topologies on a locally quasi-convex abelian group and the Mackey group problem, arXiv:2204.04389 [math.GR] (2022). \href{https://arxiv.org/abs/2204.04389}{arXiv:2204.04389}.
\bibitem{DTY2005} D. Dikranjan, M. Tkachenko, and I. Yaschenko, On transversal group topologies, Topology and its Applications 153, nos. 5--6 (2005), 786--817. \href{https://doi.org/10.1016/j.topol.2005.01.011}{doi:10.1016/j.topol.2005.01.011}.

\bibitem{Barrera2014} D. de la Barrera Mayoral, $\Q$ is not a Mackey group, Topology and its Applications 178 (2014), 265--275. \href{https://doi.org/10.1016/j.topol.2014.10.004}{doi:10.1016/j.topol.2014.10.004}.


\bibitem{AD2021} L. Außenhofer and D. Dikranjan, Mackey groups and Mackey topologies, Dissertationes Mathematicae 567 (2021), 1--141. \href{https://doi.org/10.4064/dm835-7-2021}{doi:10.4064/dm835-7-2021}.

\bibitem{AB2012} L. Außenhofer and D. de la Barrera Mayoral, Linear topologies on $\Z$ are not Mackey topologies, Journal of Pure and Applied Algebra 216:6 (2012), 1340--1347. \href{https://doi.org/10.1016/j.jpaa.2012.01.005}{doi:10.1016/j.jpaa.2012.01.005}. Theorem 4.4 in the checked preprint: \href{https://arxiv.org/abs/1107.2661}{arXiv:1107.2661}.

\bibitem{Gabriyelyan2018} S. Gabriyelyan, A locally quasi-convex abelian group without a Mackey group topology, Proceedings of the American Mathematical Society 146 (2018), 3627--3632. \href{https://doi.org/10.1090/proc/14020}{doi:10.1090/proc/14020}. \href{https://arxiv.org/abs/1708.06107}{arXiv:1708.06107}.

\end{thebibliography}
\end{document}